\documentclass[a4paper,11pt]{amsart}

\usepackage{amsfonts,amsthm,amssymb,amsmath,amscd}

\title[Pentagonal curves]{On the Family of Pentagonal curves of 
genus 6 and associated modular forms on the Ball}
\author{Kenji Koike}
\address{Department of Mathematics \& Informatics, Faculty of Science,
CHIBA UNIVERSITY, 1-33 Yayoi-cho, Inage-ku, Chiba 263-8522,Japan}
\date{}
\email{mkoike@math.s.chiba-u.ac.jp}
\numberwithin{equation}{section}
\newtheorem{lem}{Lemma}[section]         
\newtheorem{defin}{Definition}[section] 
\newtheorem{prop}{Proposition}[section] 
\newtheorem{cor}{Corollary}[section] 
\newtheorem{rem}{Remark}[section] 
\newtheorem{tm}{Theorem}[section] 

\newcommand{\Z}{\mathbb{Z}}
\newcommand{\C}{\mathbb{C}}
\newcommand{\pr}{\mathbb{P}^1}
\newcommand{\x}{X(2,5)}
\newcommand{\xc}{X^{\circ}(2,5)}
\newcommand{\B}{\mathbb{B}_2^A}
\newcommand{\Bc}{\mathbb{B}^{\circ}_2}
\newcommand{\z}{\zeta}
\newcommand{\s}{\mathfrak{S}}
\newcommand{\per}{\mathrm{S}_5}

\begin{document}
 
\maketitle

\begin{abstract}
In this article we study the inverse of the period map for
the family  $\mathcal{F}$ of complex algebraic curves of genus 6 equipped with
an automorphism of order 5. This is a family with 2 parameters, and is 
fibred over a certain type of Del Pezzo surace.
The period satisfies the hypergeometric differential equation for Appell's 
$F_1(\frac{3}{5},\frac{3}{5},\frac{2}{5},\frac{6}{5})$ of two variables after
a certain normalization of the variable parameter. \\
\indent
This differential equation and the family $\mathcal{F}$ are studied by 
G. Shimura (1964), T. Terada (1983, 1985), P. Deligne - G.D. Mostow
 (1986) and T. Yamazaki- M. Yoshida(1984). Recently  M. Yoshida 
presented a new approch using the concept of configration space.
Based on their results we show the representation of the inverse of the
period map in terms of Riemann theta constants. This is the first variant
of the work of H. Shiga (1981) and K. Matsumoto (1989, 2000) to the 
co-compact case.
\end{abstract}


\setcounter{section}{-1}
\section{Introduction}
Let $\mathcal{F}$ be the family of algebraic curves given by
\[
C(\lambda ):\  w^5=\prod_{i=1}^5(z-\lambda _i),
\]
here the parameter 
$\lambda = (\lambda_1, \lambda_2,\lambda_3, \lambda_4,\lambda_5)$ 
lives on the domain $(\pr)^5 - \Delta$, where¡¡
\[
\Delta = \{(\lambda_i) \in (\pr)^5 \ : \ \lambda_i = \lambda_j 
 \ \text{for some} \ i \ne j\}.
\]
By putting $(\lambda _1, \lambda _2, \lambda _3)=(0,1,\infty )$ we can
normalize $C(\lambda )$ in the form
\[
C'(x,y) : w^5=z(z-1)(z-x)(z-y)
\]
where the parameter $(x,y)$ lives in
\[
\Lambda = \{(x,y) \in \C^2 \ : \ xy(x-1)(y-1)(x-y) \ne 0\}
\]
The period of $C'(x,y)$
\[
\eta (x,y)=\int_{\gamma }\frac{dz}{w^2}
\]
satisfies the system of differential equation
\begin{equation} \label{appell}
\begin{split}
x(1-x) \frac{\partial^2 u}{\partial x^2} + 
y(1-x) \frac{\partial^2 u}{\partial x \partial y} + 
(\frac{6}{5} - \frac{11}{5}x) \frac{\partial u}{\partial x} - 
\frac{3}{5} y \frac{\partial u}{\partial y} - \frac{9}{5} u &= 0 \\
y(1-y) \frac{\partial^2 u}{\partial y^2} + 
x(1-y) \frac{\partial^2 u}{\partial x \partial y} + 
(\frac{6}{5} - 2 y) \frac{\partial u}{\partial y} - 
\frac{2}{5} x \frac{\partial u}{\partial x} - \frac{6}{5} u &= 0
\end{split}
\end{equation}
It is the hypergeometric differential equation for the Appell's 
$F_1(\frac{3}{5},\frac{3}{5},\frac{2}{5},\frac{6}{5};x,y)$. 
The dimension of the solution space is equal to 3. 
If it holds $\lambda' = g \circ \lambda$ for a certain projective 
transformation $g \in \mathrm{PGL}_2(\C)$, then we have the biholomorphic 
equivalence $C(\lambda) \cong C(\lambda')$. So we consider the quotient
space  
\[
\xc = ((\pr)^5 - \Delta) / \mathrm{PGL}_2(\C). 
\]
as the parameter space for $\mathcal{F}$, that is biholomorphically 
equivalent with $\Lambda $.
\par
According to the work of T. Terada (\cite{Te}),  
P. Deligne - G.D. Mostow (\cite{DM}) and T. Yamazaki - M. Yoshida(\cite{YY})
we have the following properties :
\begin{enumerate}
\item Let $\{ \eta_1, \eta_2, \eta_3 \}$ be the basis of the 
solutions of (\ref{appell}). The image of the Schwarz map 
$(x,y)\mapsto [\eta_1(x,y):\eta_2(x,y):\eta_3(x,y)]\in {\mathbb P}^2$ 
is an open dense subset of a 2-dimensional ball ${\mathbb B}_2$.
\item  The monodromy group for (\ref{appell}) is characterized as a certain 
      congruence subgroup of the Picard modular group for 
      $k={\mathbb Q}(e^{2\pi i/5})$.
\item  Let $\per$  be the symmetric group of permutations of 
$\{ \lambda_1,\lambda_2,\lambda_3,\lambda_4,\lambda_5 \}$, 
it has a natural action on $\xc$. There is a compactification $\x$ of
$\xc$ so that we have $\per \subset \mathrm{Aut}(\x)$. 
Yoshida showed $\x$ is a Del Pezzo surface of degree 5.
\item We obtain a single valued modular map on ${\mathbb B}_2$ as the 
inverse of the Schwarz map.
\end{enumerate}
By the so-called Picard principle we can reduce the period map for
 $\mathcal{F}$ to the Schwarz map for (\ref{appell}).
So we proceed our study by the following steps.
\par
In first 4 sections we make up the explicit realization of the above 
properties (1) - (4):
\par
\par
Section 1. We describe the parameter space $\xc$ for $\mathcal{F}$
and its compactification $\x$. We list up certain divisors those become to be
essential in our study.
\par
Section2. We construct the period map for ${\mathcal F}$. 
And we show how it reduces to the map $\Phi : \xc \rightarrow {\mathbb B_2}$.
\par
Section 3. We list up the generator system of the monodromy group for
$\Phi$ in terms of the unitary reflections.
\par
Section 4. We observe the degeneration of the Schwarz map for (\ref{appell}).
\par
Section 5 is the main part of the article.
There we study the 0 values of the Riemann theta functions of genus 6 with 
the characteristic $(a,b)\in (\frac{1}{10}\Z)^6\times
(\frac{1}{10}\Z)^6$ (Theorem \ref{Gaussinverse}).
These are considered to be a certain kind of automorphic form on 
${\mathbb B}_2$. Many of the above theta constants identically vanish 
on ${\mathbb B}_2$. At first it is proved that there are only 25 among
 them those are
invariant under the action of the monodromy group. We show they are not
identically zero on $\B$.
Then we determine the vanishing locus on ${\mathbb B}_2$
of every theta constant in question.
\par
In Section 6 we state the main theorem (Theorem \ref{maintheorem}), 
that is the representation of $\Phi ^{-1}$ via the theta constants.
As the direct consequence we show the representation of the inverse 
Schwarz map for the Gauss hypergeometric differential equation 
$E_{2,1}(\frac{1}{5}, \frac{2}{5}, \frac{4}{5})$. In this case we have the
arithmetic triangle group of co-compact type $\Delta (5,5,5)$ as the
monodromy group, and it is the
 case mentioned by Shimura \cite{Sm}. As the another application we show the
explicit generator system for the graded ring of the automorphic forms 
with respect to the unitary group 
$U(2,1;{\mathcal O}_k)$ over ${\mathcal O}_k$ (Theorem \ref{tmmodring}).
\\
\section{The configuration space $\x$} \label{sectionconfig}
Here we summarize the fundamental facts of $\x$.
For precise arguments, see \cite[Chapter V]{Y2}. Let $[a:b]$ be a point
on $\pr$, and let $\lambda = b/a$ be its representative on $\C
\cup \{\infty\}$. Always we use the notation $\lambda_i \in \pr$ in this
sense. Let us consider ordered distinct five
points on $\pr$:
\[
  \lambda = (\lambda_1, \lambda_2,\lambda_3, \lambda_4,\lambda_5) 
\in (\pr)^5 - \Delta
\]
where, $\Delta$ is degenerate locus:
\[
 \Delta = \{(\lambda_i) \in (\pr)^5 \ : \ \lambda_i = \lambda_j 
 \ \text{for some} \ i \ne j\}.
\]
A projective transformation $g \in \mathrm{PGL}_2(\C)$ acts on $(\pr)^5$ as
\[
  g \cdot (\lambda_1, \cdots, \lambda_5) = (g(\lambda_1), \cdots,
  g(\lambda_5)).
\]
The configuration space $\xc$ is defined by the quotient space
\[
\xc = ((\pr)^5 - \Delta) / \mathrm{PGL}_2(\C). 
\]
It has a good compactification
\[
\x = \overline{\xc} = ((\pr)^5 - \Delta') / \mathrm{PGL}_2(\C)
\]
where
\[
\Delta' =  \{(\lambda_i) \in (\pr)^5 \ : \ \lambda_i = \lambda_j = \lambda_k
 \ \text{for some} \ i \ne j \ne k \ne i\}.
\]
There exist ten lines on $\x$ of the form
\[
 L(ij) = \{\text{the orbit of the form $\lambda_i = \lambda_j$} \} / 
 \mathrm{PGL}_2(\C) \cong \pr.
\]
Notice that $L(ij) \cap L(jk) = \phi \ (i \ne j \ne k \ne i)$ 
by the definition, 
and the degenerate locus $\x - \xc$ is just the union of these
ten lines.  $\x$ is isomorphic to the blow-up of $\mathbb{P}^2$ at
four points. We can see the blow down 
$\pi : \x \rightarrow \mathbb{P}^2$ by the following way. 
Let us specialize $\lambda_4 = 0, \ \lambda_5 = \infty$ and 
regard $[\lambda_1:\lambda_2:\lambda_3]$ as a 
point in $\mathbb{P}^2$, then we obtain the following correspondence;
\begin{align*}
 P_1 = [1:0:0] = \pi( L(15)), \quad 
 P_2 = [0:1:0] = \pi( L(25)), \\
 P_3 = [0:0:1] = \pi( L(35)), \quad 
 P_4 = [1:1:1] = \pi( L(45)), 
\end{align*}
and 
\[
\pi (\xc) = \{[\lambda_1:\lambda_2:\lambda_3] \in \mathbb{P}^2 \ : \
\lambda_i \ne \lambda_j \ (i \ne j), \quad i,j = 1,2,3,4\}.
\]
\\
\indent For five distinct numbers $i,j,k,l,m$ in $\{1,2,3,4,5\}$, We define 
a divisor $D(ijklm)$ on $\x$ by
\[
 D(ijklm) = L(ij) + L(jk) + L(kl) + L(lm) + L(mi).
\]
Such a divisor is understood as a ``juzu sequence''(see \cite{Y2}). 
A 5-juzu sequence $(ijklm)$ is the pentagon with vertices $i,j,k,l,m$ in 
this cyclic order. 
The divisor $D(ijklm)$ is given by the edges of this pentagon.
\begin{figure}[htbp] \begin{center}
\setlength{\unitlength}{1mm}
\begin{picture}(35,30)
\put(10,5){\line(1,0){15}} \put(29,8){\line(1,4){2.5}}
\put(20,30){\line(2,-1){10}} \put(5,25){\line(2,1){10}} 
\put(4,18){\line(1,-4){2.5}}
\put(17,30){i} \put(3,20){j} \put(7,4){k} \put(27,4){l} \put(30,20){m}
\end{picture} \qquad
\begin{picture}(35,30)
\put(10,5){\line(1,0){15}} \put(29,8){\line(1,4){2.5}}
\put(20,30){\line(2,-1){10}} \put(5,25){\line(2,1){10}} 
\put(4,18){\line(1,-4){2.5}}
\put(17,30){i} \put(3,20){m} \put(7,4){l} \put(27,4){k} \put(32,20){j}
\end{picture}
\end{center} 
\caption{} \label{Juzu}
\end{figure}
We identify $(ijklm)$ and $(imlkj)$ since $L(ij) = L(ji)$.
There are twelve different divisors of this form. Let $H$ be a line on 
$\mathbb{P}^2$. As easily shown, $D(ijklm)$ are linearly equivalent to
the divisor
\[
  3 \pi^*H - L(15) - L(25) - L(35) - L(45).
\]
By the general theory of Del Pezzo surfaces 
(for example, see \cite[Chapter 5]{F}), 
this is anti-canonical class and very ample. 
In fact, we have the following proposition by direct calculations.
\vskip3mm
\begin{prop} \label{propJ}
Set
\begin{align*}
J(ijklm)(\lambda) = (\lambda_i - \lambda_j)(\lambda_j - \lambda_k)
(\lambda_k - \lambda_l)(\lambda_l - \lambda_m)(\lambda_m - \lambda_i)
\end{align*}
for twelve $(ijklm)$. Then the map
\begin{align*}
J \ : \ \x \longrightarrow \mathbb{P}^{11}, \quad 
J(\lambda) = [\cdots : J(ijklm)(\lambda) : \cdots]
\end{align*}
is an embedding.
\end{prop}
\begin{rem}
It is necessary to make precise the above notation for $J(ijklm)$. By
 using the homogeneous coordinate $[a_i:b_i]$ for $\lambda_i \in \pr$,
 we set $d(ij) = a_j b_i - a_i b_j$. So $\lambda_i - \lambda_j$ stands
 for $d(ij)$. The ratio $[J(ijklmn):J(i'j'k'l'm')]$ defines a rational
 function on $X(2,5)$.
\end{rem}
We shall give the correspondence between these divisors and theta 
functions in a later section.
\\
\\
\section{The family of pentagonal curves and the periods} 
\label{sectionperiod}
Let us consider the algebraic curve
\begin{align*}
 C_{\lambda} \ : \ y^5 = 
 (x -\lambda_1)(x -\lambda_2)(x -\lambda_3)(x -\lambda_4)(x -\lambda_5),
 \cr \lambda = (\lambda_1,\lambda_2,\lambda_3,\lambda_4,\lambda_5) \in 
(\pr)^5 - \Delta.
\end{align*}
If it holds $\lambda' = g \cdot \lambda$ for some $g \in
\mathrm{PGL}_2(\C)$, then we have a biholomorphic equivalence
$C(\lambda) \cong C(\lambda')$. So we identify $C(\lambda)$ and
$C(\lambda')$ in this case. 
Set $\mathcal{F} = \{ C_{\lambda} : \lambda \in \xc\}$.
We regard $C_{\lambda}$ as a five sheeted
cyclic covering over $\pr$ branched at $\lambda_i$ via the projection
\[
\pi \ : \ C_{\lambda} \longrightarrow \pr, \quad (x,y) \mapsto x.
\]
By the Hurwitz formula, the genus of $C_{\lambda}$ is six.
We have the following basis of $\mathrm{H}^0(C_{\lambda},\Omega^1)$:
\begin{align} \label{basiof1forms}
\varphi_1 = \frac{\mathrm{d}x}{y^2}, \quad 
\varphi_2 = \frac{\mathrm{d}x}{y^3}, \quad
\varphi_3 = \frac{x\mathrm{d}x}{y^3}, \quad 
\varphi_4 = \frac{\mathrm{d}x}{y^4}, \quad
\varphi_5 = \frac{x\mathrm{d}x}{y^4}, \quad 
\varphi_6 = \frac{x^2\mathrm{d}x}{y^4}.
\end{align}
Let $\rho$ denotes the automorphism of order five:
\[
 \rho \ : \ C_{\lambda} \longrightarrow C_{\lambda}, \quad (x,y) 
\mapsto (x,\zeta y) \qquad (\zeta = \exp(2 \pi \sqrt{-1}/5))
\]
on $C_{\lambda}$.
\begin{rem}
Throughout this article always $\z$ stands for $\exp(2 \pi \sqrt{-1}/5)$.
\end{rem}
Next, we construct a symplectic basis of
$\mathrm{H}_0(C_{\lambda},\Z)$. \\
Let $\lambda^0 = (\lambda_1^0,\lambda_2^0,\lambda_3^0,
\lambda_4^0,\lambda_5^0) \in \xc$ be a real point such that 
$\lambda_1 < \cdots <\lambda_5$ and $C_0$ be the corresponding curve.
Take a point $x_0 \in \pr$ such that $\mathrm{Im}(x_0) < 0$, and make
line segments $l_i$ connecting 
$x_0$ and $\lambda_i$. Then $\Sigma = \pr - \cup l_i$ is simply
connected and $\pi^{-1}(\Sigma)$ is isomorphic to $\Sigma \times \Z / 5 \Z$.
Here, we choose the fiber coordinates $k \in \Z / 5 \Z$ such that
$\rho(x,k) = (x,k+1)$.
Let $\alpha(i,j)$ be the oriented arc from $\lambda_i$ to $\lambda_j$ in
$\Sigma$. We obtain the following five oriented arcs $\alpha_k(i,j) \ (k =
1,\cdots,5)$ in $C_0$:
\begin{align} \label{arc}
 \alpha_k(i,j) = (\alpha(i,j), k) \subset \Sigma \times \Z / 5 \Z.
\end{align}
We define cycles $\gamma_1,\ \gamma_2,\ \gamma_3$ on $C_0$ 
(Figure \ref{Fgcycle}) using this notation;
\begin{equation} \label{gamma}
\begin{split}
 \gamma_1 &= \alpha_1(1,2) + \alpha_2(2,1), \cr
 \gamma_2 &= \alpha_1(3,4) + \alpha_2(4,3), \cr
 \gamma_3 &= \alpha_1(1,3) + \alpha_2(3,4) + \alpha_3(4,2) + \alpha_2(2,1).
\end{split} \end{equation}
\begin{figure}[htbp] \begin{center}
\setlength{\unitlength}{1mm}
\begin{picture}(110,45)
\thicklines
\put(10,10){\line(0,1){20}} \put(60,0){\line(-5,1){50}}
\put(35,10){\line(0,1){20}} \put(60,0){\line(-5,2){25}}
\put(60,0){\line(0,1){30}} \put(85,10){\line(0,1){20}}
\put(60,0){\line(5,2){25}} \put(110,10){\line(0,1){20}}
\put(60,0){\line(5,1){50}}
\put(10,30){\oval(5,10)[l]} \put(10,35){\vector(1,0){10}}
\put(20,35){\line(1,-2){5}} \put(25,25){\line(1,0){10}}
\thinlines
\put(35,30){\oval(5,10)[r]} \put(10,25){\line(1,0){10}}
\put(20,25){\line(1,2){5}} \put(25,35){\line(1,0){10}}
\thicklines
\put(60,30){\oval(5,10)[l]} \put(60,35){\vector(1,0){10}}
\put(70,35){\line(1,-2){5}} \put(75,25){\line(1,0){10}}
\thinlines
\put(85,30){\oval(5,10)[r]} \put(60,25){\line(1,0){10}}
\put(70,25){\line(1,2){5}} \put(75,35){\line(1,0){10}}
\thicklines
\put(10,30){\oval(10,30)[l]} \put(10,45){\vector(1,0){30}}
\put(40,45){\line(1,-2){15}} \put(55,15){\line(1,0){5}}
\put(85,30){\oval(10,30)[r]} \put(55,45){\line(1,0){30}}
\put(40,15){\line(1,2){15}} \put(35,15){\line(1,0){5}}
\thinlines
\put(60,15){\line(1,0){25}} \put(10,15){\line(1,0){25}}
\put(15,5){$l_1$} \put(32,7){$l_2$} \put(56,5){$l_3$}
\put(70,6){$l_4$} \put(102,5){$l_5$}
\put(10,37){$\gamma_1$} \put(60,37){$\gamma_2$} \put(45,37){$\gamma_3$}
\put(10,30){$\lambda_1$} \put(30,30){$\lambda_2$} \put(60,30){$\lambda_3$}
\put(80,30){$\lambda_4$} \put(105,30){$\lambda_5$}
\end{picture}
\end{center} 
\caption{} \label{Fgcycle}
\end{figure}

We set
\begin{equation} \label{AB}
\begin{split}
 &A_1 = \gamma_1, \quad A_2 = \gamma_2, \quad A_3 = \gamma_3, \quad 
 A_4 = \rho^2(\gamma_1), \quad A_5 = \rho^2(\gamma_2), \quad
 A_6 = \rho^4(\gamma_3), \\
 &B_1 = \rho(\gamma_1) + \rho^3(\gamma_1), \quad 
 B_2 = \rho(\gamma_2) + \rho^3(\gamma_2), \quad
 B_3 = \rho(\gamma_3) + \rho^2(\gamma_3), \\
 &B_4 = \rho^3(\gamma_1), \quad B_5 = \rho^3(\gamma_2), \quad
 B_6 = \rho(\gamma_3).
\end{split}
\end{equation}
The intersection numbers of these cycles are given by
\[
 A_i \cdot A_j = B_i \cdot B_j = 0, \quad A_i \cdot B_j = \delta_{ij}.
\]
So, $\{A_i,\ B_i \}$ is a symplectic basis of $\mathrm{H}_1(C_0,\Z)$. 
Let $\lambda$ be a point on $\xc$, and suppose an arc $r$ from
$\lambda^0$ to $\lambda$. Since the family $\mathcal{F}$ is locally 
trivial as a topological fiber space over $\xc$, by using this
trivialization along $r$, we obtain the systems $\{\alpha_k(i,j)(\lambda)\}$, 
$\{\gamma_i(\lambda)\}$ and the symplectic basis 
$\{A_i(\lambda),\ B_i(\lambda) \}$ on $C_{\lambda}$. We have the
relation (\ref{AB}) between 
$\{\gamma_i(\lambda)\}$ and $\{A_i(\lambda),\ B_i(\lambda) \}$ also. 
We note that $\{A_i(\lambda),\ B_i(\lambda) \}$ depend on the homotopy
class of $r$.
\\
Now, we consider the period matrix of $C_{\lambda}$:
\begin{align*}
\Pi(\lambda) = \Pi = (Z_1, Z_2) = \begin{pmatrix}
\int_{A_1} \varphi_1 & \cdots & \int_{A_6} \varphi_1 &
\int_{B_1} \varphi_1 & \cdots & \int_{B_6} \varphi_1 \\
\hdotsfor{6} \\
\int_{A_1} \varphi_6 & \cdots & \int_{A_6} \varphi_6 &
\int_{B_1} \varphi_6 & \cdots & \int_{B_6} \varphi_6
\end{pmatrix}.
\end{align*}
The normalized period matrix $\Omega(\lambda) = \Omega = Z_1^{-1} Z_2$ 
belongs to the Siegel upper half space of degree 6:
\[
 \mathfrak{S}_6 = \{\Omega \in \mathrm{GL}_6(\C) \ : \ {}^t \Omega = \Omega,\
 \mathrm{Im}(\Omega) \ \text{is positive definite} \}.
\]
The automorphism $\rho$ acts on $\mathrm{H}^0(C_{\lambda},\Omega^1)$ and
$\mathrm{H}_1(C_{\lambda},\Z)$. So we have the representation matrices $R
\in \mathrm{GL}_6(\C)$ and $M \in \mathrm{GL}_{12}(\Z)$ of $\rho$ 
with respect to the basis $\{\varphi_i\}$ and $\{A_i, B_i\}$, respectively.
It holds $R \Pi = \Pi M$. Put
\begin{align} \label{sigma}
 M = \begin{pmatrix} {}^tD & {}^tB \\ {}^tC & {}^tA \end{pmatrix}, \quad
 \sigma = \begin{pmatrix} A & B \\ C & D \end{pmatrix}.
\end{align}
Then the matrix $\sigma$ belongs to the symplectic group 
\[
 \mathrm{Sp}_{12}(\Z) = \{g \in \mathrm{GL}_{12}(\Z) \ : \ {}^tg J g = J
 \}, \quad J = \begin{pmatrix} 0 & I_6 \\ -I_6 & 0 \end{pmatrix}
\]
and it holds
\[\Omega = (A \Omega + B)(C \Omega + D)^{-1}.
\]
As easily shown, $\varphi_i \ (i = 1, \cdots, 6)$ is eigenvectors of 
$\rho$ and we have 
\[
 R = \begin{pmatrix}
\zeta^3 &&&&& \\ & \zeta^2 &&&0& \\ && \zeta^2 &&& \\
&&& \zeta && \\ &0&&& \zeta & \\ &&&&& \zeta
\end{pmatrix}.
\]
By the relation (\ref{AB}) of $A_i, B_i$, we have
\begin{align}\label{piri}
 \Pi = (a,b,c,R^2 a,R^2 b,R^4 c,(R + R^3)a,(R + R^3)b,(R + R^2)c,R^3 a,
 R^3 b,R c),
\end{align}
where we denote
\[a = {}^t(\int_{\gamma_1} \varphi_1, \cdots, \int_{\gamma_1}\varphi_6),\quad
  b = {}^t(\int_{\gamma_2} \varphi_1, \cdots, \int_{\gamma_2}\varphi_6),\quad
  c = {}^t(\int_{\gamma_3} \varphi_1, \cdots, \int_{\gamma_3}\varphi_6). \]
According to (\ref{AB}),
\[
 \rho (A_1) =\rho (\gamma_1) = (\rho (\gamma_1) + \rho^3 (\gamma_1)) -
\rho^3 (\gamma_1) = B_1 - B_4.
\]
By the same way, we can describe $\rho (A_2) , \cdots, \rho (B_6)$ in 
terms of $\{A_i, B_i\}$. So we can determine $M$, and obtain
\begin{equation}\label{matrixsigma}
\sigma = \left( \begin{array}{cccccccccccc}
-1&0&0&0&0&0& -1&0&0&0&0&0 \\
0&-1&0&0&0&0& 0&-1&0&0&0&0 \\
0&0&0&0&0&-1& 0&0&-1&0&0&-1 \\ 
-1&0&0&0&0&0& -1&0&0&-1&0&0 \\
0&-1&0&0&0&0& 0&-1&0&0&-1&0 \\
0&0&1&0&0&-1& 0&0&0&0&0&0 \\
1&0&0&-1&0&0& 0&0&0&0&0&0 \\
0&1&0&0&-1&0& 0&0&0&0&0&0 \\
0&0&0&0&0&1& 0&0&0&0&0&0 \\
0&0&0&1&0&0& 0&0&0&0&0&0 \\
0&0&0&0&1&0& 0&0&0&0&0&0 \\
0&0&0&0&0&0& 0&0&1&0&0&0 \\
\end{array} \right). 
\end{equation}
Put
\[
 \eta(\lambda) = [\eta_1(\lambda):\eta_2(\lambda):\eta_3(\lambda)] 
\in \mathbb{P}^2, \quad
  \eta_1(\lambda)  = \int_{\gamma_1} \varphi_1, \quad
  \eta_2(\lambda) = \int_{\gamma_2} \varphi_1, \quad
  \eta_3(\lambda) = \int_{\gamma_3} \varphi_1. 
\]
These are multi-valued analytic functions of $\lambda$.
Applying the Riemann positive condition 
\[
 (\int_{A_1} \varphi_1, \cdots, \int_{B_6} \varphi_1) \ J \
{}^t(\int_{A_1} \overline{\varphi}_1, \cdots, 
\int_{B_6} \overline{\varphi}_1) > 0
\]
for (\ref{piri}), we obtain
\[
 |\eta_1|^2 + |\eta_2|^2 + \frac{1 - \sqrt{5}}{2} |\eta_3|^2 < 0.
\]
So, $\eta = (\eta_1, \eta_2, \eta_3)$ belongs to the complex ball
\begin{align} \label{ballA}
 \B = \{\eta \in \mathbb{P}^2 \ : \ {}^t\bar{\eta}A\eta < 0\}, \quad
 A = \mathrm{diag}(1,1,\frac{1 - \sqrt{5}}{2}).
\end{align}
Next, we determine $\Omega$ explicitly.
Write $a = (a_i), \ b = (b_i)$ and $c = (c_i)$. Then, the Riemann bilinear
relation $\Pi J {}^t\Pi = 0$ induces the following equations:
\[
 c_2 = -(\zeta^2 + \zeta^3)(a_1 a_2 + b_1 b_2) / c_1, \quad
c_3 = -(\zeta^2 + \zeta^3)(a_1 a_3 + b_1 b_3) / c_1.
\]
By substituting them for $Z_1, \ Z_2$ in $\Pi$ we can proceed the calculation 
 of $\Omega = Z_1^{-1} Z_2$ (using a computer). Hence we have the 
following:
\begin{lem}\label{Period}
Let $\Delta = \eta_1^2 + \eta_2^2 - \z^3 (1 + \z) \eta_3^2$.
The period matrix $\Omega = (\Omega_{ij})$ is given by
\begin{align*}
\Omega_{11} &=& (\z^3-1)(\eta_1^2 +(1+\z^3) \eta_2^2 +\eta_3^2)/\Delta, \quad
\Omega_{44} &=& -\z^2 (\eta_1^2 + \z^2 \eta_2^2 -(1+\z) \eta_3^2)/\Delta, \\
\Omega_{22} &=& (\z^3-1)((1+\z^3) \eta_1^2 +\eta_2^2 + \eta_3^2)/\Delta, \quad
\Omega_{55} &=& -\z^2(\z^2\eta_1^2 +\eta_2^2 - (1+\z)\eta_3^2)/\Delta, \\
\Omega_{33} &=& (\z^2 -1)(\eta_1^2 +\eta_2^2 - \z^3 \eta_3^2)/\Delta, \quad
\Omega_{66} &=& -\z^3( \eta_1^2 + \eta_2^2 -(1+\z^4)\eta_3^2)/\Delta, \\
\Omega_{12} &=& (\z^3 - \z) \eta_1 \eta_2/\Delta, \quad 
\Omega_{45} &=& (\z^4 - \z^2) \eta_1 \eta_2 / \Delta, \\
\Omega_{15} &=& (\z^4 - \z) \eta_1 \eta_2 / \Delta, \quad
\Omega_{24} &=& (\z^4 -\z) \eta_1 \eta_2 / \Delta, \\
\Omega_{13} &=& (1 - \z^2) \eta_1 \eta_3 / \Delta, \quad
\Omega_{23} &=& (1 - \z^2) \eta_2 \eta_3 / \Delta, \\
\Omega_{46} &=& (\z^4 - \z) \eta_1 \eta_3 / \Delta, \quad
\Omega_{56} &=& (\z^4 - \z) \eta_2 \eta_3 / \Delta, \\
\Omega_{16} &=& (\z^3 - \z) \eta_1 \eta_3 / \Delta, \quad
\Omega_{26} &=& (\z^3 - \z) \eta_2 \eta_3 / \Delta, \\
\Omega_{34} &=& (1-\z^3) \eta_1 \eta_3 / \Delta, \quad
\Omega_{35} &=& (1-\z^3) \eta_2 \eta_3 / \Delta, \\
\Omega_{14} &=& \z^3((1 + \z) \eta_1^2 + (1 + \z^3) \eta_2^2 + \eta_3^2) / 
\Delta, \quad
\Omega_{25} &=& \z^3((1 + \z^3) \eta_1^2 + (1 + \z) \eta_2^2 + \eta_3^2) / 
\Delta, \\
\Omega_{36} &=& (\z + \z^2)(\eta_1^2 + \eta_2^2 - 
\zeta^3 (1 + \zeta^2) \eta_3^2) / \Delta .
\end{align*}
\end{lem}
\vskip3mm
Now we define our period map
\[
 \Phi \ : \ \xc \longrightarrow \B, \quad \lambda 
\mapsto [\eta_1(\lambda) : \eta_2(\lambda) : \eta_3(\lambda)],
\]
that is multi-valued analytic. The above Lemma says that the original 
period map $\lambda \mapsto \Omega(\lambda)$ factors as
\[
 \xc \longrightarrow \B \longrightarrow \s_6.
\]
Throughout this paper, we denote matrices of the form in Lemma
\ref{Period} by $\Omega(\eta)$.
\\
\section{The monodromy group and Reflections}
The multi-valuedness of $\Phi$ induces a unitary representation with
respect to $A$ in (\ref{ballA}) of the fundamental group $\pi_1(\xc)$.
We call it the monodromy group of $\Phi$. The structures of our
monodromy group is studied in \cite{YY}. Set
\[
 \Gamma = \{g \in \mathrm{GL}_3(\Z[\z]) \ : \ {}^t \bar{g} A g  = A \},
\quad 
\Gamma(1-\z) = \{g \in \Gamma \ : \ g \equiv I_3 \ \mathrm{mod} \ 1-\z
\}.
\]
The group $\Gamma$ acts on $\B$ (left action).
\begin{tm}[T. Yamazaki, M. Yoshida \cite{YY}]\label{tmYY}
(1). \  The monodromy group of the period map $\Phi$ coincides with 
$\Gamma(1-\z)$ and the quotient $\Gamma / (\pm I) \Gamma(1-\z)$ is 
isomorphic to the symmetric group $\per$. \\
(2). The quotient $\B / \Gamma(1-\z)$ is biholomorphically equivalent 
to the blow up of $\mathbb{P}^2$ at four points.
\end{tm}
\begin{rem}[see \cite{YY}]
There are ten $(-1)$--curves on $\B / \Gamma(1-\z)$, 
and $\per$ acts transitively on them.
\end{rem}
According to \cite{Te} and \cite{YY}, it is proved that $\Gamma$ and 
$\Gamma(1-\z)$ are reflection groups and the generator systems are given 
also. We expose those generator system in a form adapted for our
calculation in the later sections. \\
\indent Let us consider the reference point $\lambda^0 \in \xc$ again.
Now we define the half way monodromy transformation $g_{12}$ induced 
from the permutation of $\lambda^0_1$ and $\lambda^0_2$. 
Let us consider a continuous arc $R_{12}$ starting from $\lambda^0$:
\begin{align} \label{lambdadeformation}
 \lambda(t) = (\lambda_1(t), \lambda_2(t),
 \lambda_3^0,\lambda_4^0,\lambda_5^0), \quad (0 \leq t \leq 1)
\end{align}
such that (Figure \ref{Fghalf})
\[
\lambda_2(1) = \lambda_1^0, \quad \lambda_1(1) = \lambda_2^0, \qquad
\mathrm{Im}(\lambda_1(t)) < 0 < \mathrm{Im}(\lambda_2(t))
\quad (0 < t < 1).
\]
\begin{figure}[htbp] \begin{center}
\setlength{\unitlength}{1mm}
\begin{picture}(110,45)
\thicklines
\put(10,10){\line(0,1){20}} \put(60,0){\line(-5,1){50}}
\put(35,10){\line(0,1){20}} \put(60,0){\line(-5,2){25}}
\put(60,0){\line(0,1){30}} \put(85,10){\line(0,1){20}}
\put(60,0){\line(5,2){25}} \put(110,10){\line(0,1){20}}
\put(60,0){\line(5,1){50}}
\put(22.5,30){\oval(25,10)[bl]} \put(22.5,30){\oval(25,10)[rt]}
\put(22.5,25){\vector(1,0){2}}  \put(22.5,35){\vector(-1,0){2}}
\put(8,31){$\lambda_1$} \put(30,31){$\lambda_2$} \put(58,31){$\lambda_3$}
\put(83,31){$\lambda_4$} \put(108,31){$\lambda_5$}
\put(25,25){$\lambda_1(t)$} \put(12,35){$\lambda_2(t)$}
\end{picture}
\end{center} 
\caption{} \label{Fghalf}
\end{figure}

Let $\eta(t) = \eta(\lambda(t))$ be the corresponding periods. 
Recall the definition (\ref{gamma}). It is apparent that $\gamma_2$ and 
$\gamma_3$ are invariant after this deformation process. Describing
$\gamma_1(t)$ for any $0 \leq t \leq 1$, we get $\gamma_1(1) = -
\rho(\gamma_1(0))$. Namely,
\[
 \left( \begin{array}{c} \eta_1(1) \\ \eta_2(1) \\ \eta_3(1) 
\end{array} \right) =
\left( \begin{array}{ccc} - \z^3 &0&0 \\ 0&1&0 \\ 0&0&1 
\end{array} \right)
\left( \begin{array}{c} \eta_1(0) \\ \eta_2(0) \\ \eta_3(0) 
\end{array} \right).
\]
The matrix in the right hand side belongs to $\Gamma$, and we denote it
by $g_{12}$.
We define $g_{23}, \ g_{34}, \ g_{45}$ by the same manner. Set
\begin{align} \label{hij}
 h_{12} = (g_{12})^2, \quad h_{23} = (g_{23})^2, 
\quad h_{34} = (g_{34})^2 \\
h_{13} = (g_{23})^{-1}(g_{12})^2g_{23},\quad 
h_{14} = (g_{23}g_{34})^{-1}(g_{12})^2g_{23}g_{34}.
\end{align}
\begin{prop}[see \cite{Y1}]
The monodromy group is generated by $h_{ij}$ in (\ref{hij}).
\end{prop}
\vskip3mm
Let $\mathrm{T}_{\alpha}$ be the reflection on $\B$ with respect 
to a root $\alpha$; 
\[
 \mathrm{T}_{\alpha}(\eta) = \eta - (1 + \z^3)
   \frac{{}^t\bar{\alpha}A\eta}{{}^t\bar{\alpha}A\alpha} \alpha,
\]
and $\mathrm{R}_{\beta}$ be the reflection on $\B$ with respect 
to a root $\beta$;
\[
\mathrm{R}_{\beta}(\eta) = \eta - (1 - \z)
   \frac{{}^t\bar{\beta}A\eta}{{}^t\bar{\beta}A\beta} \beta. 
\]
Then we see that
\begin{lem}
Set
\[
\alpha_{12} = (1,0,0), \quad \alpha_{23} = (\z^3,1,-(1+\z)), 
\quad \alpha_{34} = (0,1,0),\quad \alpha_{45} = (0,1,\z^3)
\]
and set
\begin{align*}
\beta_{12} = (1,0,0), \quad \beta_{13} = (1,-1,1+\z),
\quad \beta_{14} = (1,\z^3,1+\z), \\ 
\beta_{23} = (\z^3,1,-(1+\z)), \quad  \beta_{34} = (0,1,0).
\end{align*}
Then it holds $g_{ij} = \mathrm{T}_{\alpha_{ij}}$, 
$h_{kl} = \mathrm{R}_{\beta_{ij}}$. And $g_{ij}$ is of order five,
 $h_{kl}$ is of order ten. 
\end{lem}
\begin{rem}
The group $\Gamma'$ generated by $\{g_{ij}\}$ has a representation to
 $\per$. According to Theorem \ref{tmYY}, $\Gamma$ is generated
 $\Gamma'$ and $\pm I$.
\end{rem}
\indent The deformation of the curve $C_{\lambda(t)}$ along $R_{12}$ in
(\ref{lambdadeformation}) induces a symplectic basis $\{A_i(t), B_i(t)\}$
 on it. So $\{A_i(1), B_i(1)\}$ is again a symplectic basis on 
$C_{\lambda^0}$. Hence we obtain a symplectic transformation
\[
 {}^t(B_1(1), \cdots, B_6(1),A_1(1) \cdots,A_6(1)) = \hat{g}_{12}
{}^t(B_1(0), \cdots, B_6(0),A_1(0) \cdots,A_6(0)).
\]
For $\hat{g}_{12} = \begin{pmatrix} A&B\\C&D\end{pmatrix}$, we have
\[
  \Omega(\eta_1) = (A \Omega(\eta_0) + B)(C \Omega(\eta_0) + D)^{-1}.
\]
Recall $R_{12}$ induces the change of cycles
\[
 (\gamma_1, \gamma_2, \gamma_3) \longrightarrow 
(-\rho(\gamma_1), \gamma_2, \gamma_3).
\]
Together with (\ref{AB}), we obtain:
\begin{equation} \label{matrixg12}
\hat{g}_{12} = \left( \begin{array}{cccccccccccc}
1&0&0&0&0&0& 1&0&0&0&0&0 \\
0&1&0&0&0&0& 0&0&0&0&0&0 \\
0&0&1&0&0&0& 0&0&0&0&0&0 \\ 
1&0&0&0&0&0& 1&0&0&1&0&0 \\
0&0&0&0&1&0& 0&0&0&0&0&0 \\
0&0&0&0&0&1& 0&0&0&0&0&0 \\
-1&0&0&1&0&0& 0&0&0&0&0&0 \\
0&0&0&0&0&0& 0&1&0&0&0&0 \\
0&0&0&0&0&0& 0&0&1&0&0&0 \\
0&0&0&-1&0&0& 0&0&0&0&0&0 \\
0&0&0&0&0&0& 0&0&0&0&1&0 \\
0&0&0&0&0&0& 0&0&0&0&0&1 
\end{array} \right). 
\end{equation}

By same consideration, we obtain following;
\begin{equation}\label{matrixg23}
\hat{g}_{23} = \left( \begin{array}{cccccccccccc}
1&1&0&1&-1&1& 2&-1&-1&1&-1&1 \\
-1&1&-1&-1&1&0& -2&2&0&-1&1&-2 \\
1&-1&2&0&0&-1& 0&-1&1&0&0&1 \\
0&1&-1&1&0&1& 1&0&-1&1&0&0 \\
0&0&0&-1&1&-1& -1&1&0&-1&1&-1 \\
1&0&1&1&-1&1& 2&-2&0&1&-1&2 \\
-1&0&-1&0&1&1& 0&1&0&0&0&-1 \\
1&-1&1&0&0&-1& 0&0&1&0&0&1 \\
-1&1&-2&-1&1& 1&-1&2&0&0&1&-2 \\
0&0&0&-1&0&-1 &-1&1&0&0&1&-1 \\
0&0&1&1&-1&0& 1&-1&0&0&0&1 \\
1&-1&2&0&-1&-2 &0&-1&1&-1&0&2
\end{array} \right), 
\end{equation}
\begin{equation}\label{matrixg34}
\hat{g}_{34} = \left( \begin{array}{cccccccccccc}
1&0&0&0&0&0& 0&0&0&0&0&0 \\
0&1&0&0&0&0& 0&1&0&0&0&0 \\
0&0&1&0&0&0& 0&0&0&0&0&0 \\ 
0&0&0&1&0&0& 0&0&0&0&0&0 \\
0&1&0&0&0&0& 0&1&0&0&1&0 \\
0&0&0&0&0&1& 0&0&0&0&0&0 \\
0&0&0&0&0&0& 1&0&0&0&0&0 \\
0&-1&0&0&1&0& 0&1&0&0&0&0 \\
0&0&0&0&0&0& 0&0&1&0&0&0 \\
0&0&0&0&0&0& 0&0&0&1&0&0 \\
0&0&0&0&-1&0& 0&1&0&0&0&0 \\
0&0&0&0&0&0& 0&0&0&0&0&1 
\end{array} \right), 
\end{equation}
\begin{equation}\label{matrixg45}
\hat{g}_{45} = \left( \begin{array}{cccccccccccc}
1&0&0&0&0&0& 0&0&0&0&0&0 \\
0&1&-1&0&1&0& 0&2&0&0&1&-1 \\
0&0&1&0&-1&0& 0&0&1&0&0&1 \\ 
0&0&0&1&0&0& 0&0&0&0&0&0 \\
0&0&0&0&1&0& 0&1&1&0&1&0 \\
0&0&0&0&-1&1& 0&-1&0&0&-1&1 \\
0&0&0&0&0&0& 1&0&0&0&0&0 \\
0&-1&1&0&0&-1& 0&0&0&0&0&0 \\
0&0&-1&0&1&0& 0&1&0&0&0&-1 \\
0&0&0&0&0&0& 0&0&0&1&0&0 \\
0&0&0&0&-1&0& 0&-1&0&0&0&1 \\
0&-1&1&0&0&-1& 0&-1&0&0&0&1 
\end{array} \right). 
\end{equation}
\\
\section{Degenerate loci} \label{sectiondl}
According to Theorem \ref{tmYY}, P. Deligne - G.D. Mostow(\cite{DM}) 
and T. Terada(\cite{Te}) the period map $\Phi$ induces the biholomorphic 
equivalence
\[
 \tilde{\Phi} \ : \ \xc \stackrel{\sim}{\longrightarrow}
 \Bc / \Gamma(1 - \z),
\]
where $\Bc = \mathrm{Im} \Phi$. Moreover we have the unique extension
\[
 \tilde{\Phi} \ : \ \x \stackrel{\sim}{\longrightarrow}
 \B / \Gamma(1 - \z),
\]
and $\cup L(ij) = \x - \xc$ corresponds to $(\B - \Bc)/ 
\Gamma(1 - \z)$. Let $\pi$ be the projection 
$\B \rightarrow \B / \Gamma(1-\z)$, and let $\ell(ij)$ denote
$\pi^{-1}(\tilde{\Phi}(L(ij)))$. \\
\indent Now we consider a degenerate curve
\[
y^5 = (x-\lambda_1)^2(x-\lambda_3)(x-\lambda_4)(x-\lambda_5)
\]
with $(\lambda_1,\lambda_1,\lambda_3,\lambda_4,\lambda_5) \in L(12)$,
and putting $\lambda' = (\lambda_1,\lambda_3,\lambda_4,\lambda_5)$ we
denote it by $C_{\lambda'}$. Let $\tilde{C}_{\lambda'}$ denote the
non-singular model of $C_{\lambda'}$. It is a curve of genus 4. Set 
$\mathcal{F}_{12}$ be the totality of $\tilde{C}_{\lambda'}$. For the
parameter $(\lambda^0)' =
(\lambda_1^0,\lambda_3^0,\lambda_4^0,\lambda_5^0)$ the cycle $\gamma_1$
vanishes on $\tilde{C}_{(\lambda^0)'}$, but $\gamma_2$ and $\gamma_3$ are
still alive. So we can define $A_i, B_i \ (i = 2,3,5,6)$ on
$\tilde{C}_{\lambda'}$ by the same argument as for $C_{\lambda}$.
Hence we obtain a basis $\{A_i, B_i \} \ (i = 2,3,5,6)$ of
$\mathrm{H}_1(\tilde{C}_{\lambda'}, \Z)$. By putting 
$\lambda' = (0,1,t,\infty)$ the period 
\begin{align}
\int_{\gamma}
 x^{-\frac{4}{5}} (x-1)^{-\frac{2}{5}} (x-t)^{-\frac{2}{5}}dx,
\qquad (\gamma \in \mathrm{H}_1(\tilde{C}_{\lambda'}, \Z))
\end{align}
on $\tilde{C}_{\lambda'}$ gives a
solution for the Gauss hypergeometric differential equation 
$E_{2,1}(\frac{1}{5}, \frac{2}{5}, \frac{4}{5})$: 
\begin{align} \label{Gauss(5,5,5)}
t(1 -t) \frac{\mathrm{d}^2 u}{\mathrm{d} t^2} + 
(\frac{4}{5} - \frac{8}{5} t) \frac{\mathrm{d} u}{\mathrm{d} t}
 - \frac{2}{5} u = 0
\end{align}
The corresponding monodromy group is the triangle group $\Delta(5,5,5)$
(see \cite{Ta}, \cite{Y1}, \cite[p.138]{Y2}). Set 
\[
 \mathbb{B}_1 = \{\eta \in \B \ : \ \eta_1 = 0\},
\]
it is the mirror of the reflection $g_{12}$. By using the system
$\{\gamma_2, \gamma_3 \}$ we define a multi-valued map
\[
 \Phi_{12} \ : \ L(12) \longrightarrow \mathbb{B}_1, 
\quad \lambda \mapsto [0: \eta_2(\lambda) : \eta_3(\lambda)].
\]
It induces the restriction $\tilde{\Phi}|_{L(12)}$. By the same manner 
we obtain that $\tilde{\Phi}|_{L(ij)}$ is the mirror of the 
reflection $g_{ij}$. 
Suppose $\lambda \in L(12)$ and set $\eta = \eta(\lambda) =
\Phi_{12}(\lambda)$ . By putting $\eta_1 = 0$ in Lemma \ref{Period}, we
see that 
\begin{align} \label{degeperiod}
\Omega(\eta) =  \begin{pmatrix} \Omega_{11} &\Omega_{14}\\ 
\Omega_{41} &\Omega_{44} \end{pmatrix} \oplus \Omega'(\eta), \quad
\begin{pmatrix} \Omega_{11} &\Omega_{14}\\ 
\Omega_{41} &\Omega_{44} \end{pmatrix} = \tau_0 =
\begin{pmatrix} \z - 1 &\z^2 +\z^3 \\ 
\z^2 +\z^3 & -\z^4 \end{pmatrix}
\end{align}
with  a certain element $\Omega'(\eta) \in \s_4$. Moreover, in case
$\eta_0 = [0:0:1] \in \ell(12) \cap \ell(34)$ we have
\begin{align} \label{eta_0}
\Omega(\eta_0) = \begin{pmatrix} \Omega_{11} &\Omega_{14}\\ 
\Omega_{41} &\Omega_{44} \end{pmatrix} \oplus 
\begin{pmatrix} \Omega_{22} &\Omega_{25}\\ 
\Omega_{52} &\Omega_{55} \end{pmatrix} \oplus 
\begin{pmatrix} \Omega_{33} &\Omega_{36}\\ 
\Omega_{63} &\Omega_{66} \end{pmatrix}
= \tau_0 \oplus \tau_0 \oplus \tau_0
\end{align}
We use the above matrix to numerical evaluation of theta functions 
in later section.
\\
\\
\section{Theta functions}
\subsection{Invariant theta characteristics}
We recall basic facts on the Riemann theta functions. For a
characteristic $(a,b) \in (\mathbb{R}^g)^2$, the theta function
$\Theta_{(a,b)}(z,\Omega)$ on $\C^g \times \s_g$ is defined by the series
\[
  \Theta_{(a,b)}(z,\Omega) = \sum_{n \in \Z^g}\exp[\pi \sqrt{-1} 
{}^t(n + a)\Omega(n + a) + 2 \pi \sqrt{-1}{}^t(n + a)(z + b)].
\]
These functions satisfy the following period relations
\begin{align} \label{translate1}
\Theta_{(a,b)}(z + m,\Omega) &= \exp(2 \pi \sqrt{-1}{}^tma) 
\Theta_{(a,b)}(z,\Omega),
\end{align}
\begin{align} \label{translate2}
\Theta_{(a,b)}(z + \Omega m,\Omega) &= 
\exp(-\pi \sqrt{-1}{}^tm \Omega m -2 \pi \sqrt{-1}{}^tm(z + b))
\Theta_{(a,b)}(z,\Omega)
\end{align}
for $m \in \Z^g$. For the characteristics, we have 
\begin{align}\label{chformula}
\Theta_{(a + n,b + m)}(z,\Omega) = \exp(2 \pi \sqrt{-1}{}^tam) 
\Theta_{(a,b)}(z,\Omega),
\end{align}
for $n, \ m \in \Z^g$ and 
\begin{align}\label{chformula2}
\Theta_{(-a,-b)}(z,\Omega) = \Theta_{(a,b)}(z,\Omega).
\end{align}
Theta constants $\Theta_{(a,b)}(\Omega) = \Theta_{(a,b)}(0,\Omega)$
satisfy following transformation formula (see \cite[p176]{Ig}) 
as function on $\s_g$. For $g = \begin{pmatrix} A&B\\C&D\end{pmatrix}
 \in \mathrm{SP}_{2g}(\Z)$, set
\begin{align}
g\Omega &= (A \Omega + B)(C \Omega + D)^{-1} \\ \label{gab}
g(a,b) &= (Da - Cb,\ -Ba + Ab) + \frac{1}{2}( (C {}^tD)_0,\ (A {}^tB)_0) \\
\phi_{(a,b)}(g) &= -\frac{1}{2}(^ta^tDBa - 2^ta^tBCb + {}^tb^tCAb) + 
\frac{1}{2}(^ta^tD - ^tb^tC)(A^tB)_0
\end{align}
where $(A)_0$ stands for the diagonal vector of a matrix $A$.
Then we have
\begin{align}\label{formula}
\Theta_{g(a,b)}(g\Omega) = \kappa(g)\exp(2 \pi \sqrt{-1}\phi_{a,b}(g))
\mathrm{det}(C\Omega + D)^{\frac{1}{2}}\Theta_{(a,b)}(\Omega)
\end{align}
where, $\kappa(g)$ is a certain 8-th root of 1 depending only on $g$.
\begin{rem} \label{remtrans}
By the definition, we have
\[
\Theta_{(a,b)}(z,\Omega) = \exp(\pi \sqrt{-1} {}^ta\Omega a + 2\pi \sqrt{-1} 
{}^ta(z + b)) \Theta_{(0,0)}(z + \Omega a + b, \ \Omega), 
\]
so we often identify a characteristic $(a,b) \in (\mathbb{R}^g)^2$ with 
$\Omega a + b \in \C^g$. For 
$\begin{pmatrix} A&B \\ C&D \end{pmatrix} \in \mathrm{Sp}_{2g}(\Z)$, we
 have 
\[
 \Omega' (Da - Cb) + (-Ba + Ab) = {}^t(C \Omega + D)^{-1}(\Omega a + b),
\]
where $\Omega' = (A \Omega + B)(C \Omega + D)^{-1}$.
\end{rem}
For these formulas, see \cite{Ig} and \cite{Mum}.
\\
\\
\indent Henceforth we suppose the characteristics $(a,b)$ satisfy
$a,b \in (\frac{1}{10}\Z)^6$.
\begin{lem} \label{fixedlemma}
Let $\sigma$ be the matrix in (\ref{sigma}) and write 
$a = (a_i), \ b = (b_i)$. \\
1. \ We have
\[
 \sigma (a,b) \equiv (a,b) \ \mathrm{mod} \ \Z
\]
if and only if
\begin{align*}
5 a_1 &\equiv \frac{1}{2}, \quad a_4 \equiv a_1, \quad b_1 \equiv -2 a_1, 
\quad b_4 \equiv - a_1 \\
5 a_2 &\equiv \frac{1}{2}, \quad a_5 \equiv a_2, \quad b_2 \equiv -2 a_2, 
\quad b_5 \equiv - a_2 \qquad \mathrm{mod} \ \Z.\\
5 a_3 &\equiv \frac{1}{2}, \quad a_6 \equiv a_3, \quad b_3 \equiv -2 a_3, 
\quad b_6 \equiv - a_3 
\end{align*}
(2) \ Let $(a,b)$ be the characteristic with the above condition. Then we 
 have
\[
\hat{g} (a,b) \equiv (a,b) \ \mathrm{mod} \ \Z \quad \text{for all} \ 
g \in \Gamma(1-\z).
\]
\end{lem}
\begin{proof}
(1) \ Using the exact form (\ref{matrixsigma}) we can describe 
$\sigma (a,b)$. Then we deduce the assertion. \\
(2) \ The transformation $g(a,b)$ in (\ref{gab}) define a group action 
of the symplectic group on $(\mathbb{R} / \Z)^{2g}$ (see \cite{Ig}). 
We can check that the equality for every member of the generator system
 $\{h_{ij} \}$ of $\Gamma(1-\z)$.
\end{proof}
\vskip3mm
\begin{defin} \label{defchara}
Let $(a,b)$ be the characteristic satisfying the condition Lemma
 \ref{fixedlemma} (1). Then we can put
\begin{align}\label{charatype}
 a = \frac{1}{10}{}^t(a_1, a_2, a_3,a_1, a_2, a_3), \
 b = \frac{1}{10}{}^t(-2a_1, -2a_2, -2a_3,-a_1, -a_2, -a_3). 
\end{align}
Let $(a_1, a_2, a_3)$ denote this characteristic. We call 
$(a,b) = (a_1, a_2, a_3)$ ``$\sigma$--invariant'' if $a_1, a_2, a_3$ are 
 odd integers. For a characteristic of this type,
we denote the zero locus of $\Theta_{(a_1,a_2,a_3)}$ on $\B$ by
$\vartheta(a_1,a_2,a_3)$;
\[
 \vartheta(a_1,a_2,a_3) = \{\eta \in \B \ : \ 
\Theta_{(a_1,a_2,a_3)}(\Omega(\eta)) = 0\}.
\]
\end{defin}
\vskip3mm
\begin{rem}
By the transformation formula (\ref{formula}) and Lemma
 \ref{fixedlemma}, we see that
\[
\Theta_{(a_1,a_2,a_3)}(g \Omega(\eta)) = 
(\text{a unit function}) \times \Theta_{(a_1,a_2,a_3)}(\Omega(\eta))
\]
for a invariant characteristic $(a_1,a_2,a_3)$ and 
$g \in \Gamma(1-\z)$. Hence if we have $\eta \in \vartheta(a_1,a_2,a_3)$, 
then $\Gamma(1-\z)$--orbit of $\eta$ contained in  $\vartheta(a_1,a_2,a_3)$.
\end{rem}
\begin{lem} \label{vanish100}
Let $(a_1, a_2, a_3)$ be a $\sigma$--invariant characteristic. 
If $2 a_1^2 + 2 a_2^2 + a_3^2 \notin 5 \Z$, then 
$\vartheta(a_1,a_2,a_3) = \B$. Namely, $\Theta_{(a_1,a_2,a_3)}$ vanishes
 on $\B$. 
\end{lem}
\begin{proof}
We apply the transformation formula (\ref{formula}) for $g = \sigma^4$.
\\ For it we proceed the preparatory calculations. At first, get the
 explicit form of $g = \sigma^4 = \begin{pmatrix} A&B\\C&D\end{pmatrix}$ 
 by using (\ref{matrixsigma}). So we obtain
\[
 \phi_{(a_1,a_2,a_3)}(\sigma^4) =  \frac{1}{40}(2 a_1^2 + 2 a_2^2 + a_3^2).
\]
Using the explicit form of $\Omega(\eta)$ in Lemma \ref{Period}, we get 
\[
 \det(C\Omega(\eta) + D) = 1
\]
for all $\eta \in \B$ by a computer and calculation. By
 (\ref{chformula}), we may put
\[
\Theta_{\sigma^4(a_1,a_2,a_3)}(\Omega) = \exp[2 \pi \sqrt{-1} {}^t a m] 
\Theta_{(a_1,a_2,a_3)}(\Omega) 
\]
for a certain $m \in \Z^6$. Returning to the explicit form of 
$\sigma^4(a_1,a_2,a_3)$ we should get $m$. We check that 
$\exp[2 \pi \sqrt{-1} {}^t a m] = 1$ by a computer aided calculation.
Hence we have
\[
\Theta_{(a_1,a_2,a_3)}(\Omega(\eta)) = \kappa(\sigma^4) 
\exp[\frac{1}{20}\pi \sqrt{-1}(2 a_1^2 + 2 a_2^2 + a_3^2)]
\Theta_{(a_1,a_2,a_3)}(\Omega(\eta))
\]
for all $\eta \in \B$. This implies our assertion since
 $\kappa(\sigma^4)$ is an 8-th root of 1.
\end{proof}
\vskip 3mm
We consider odd integers $a_1,a_2,a_3$ modulo $10 \Z$.
There exist 25 representatives of the $\sigma$--invariant characteristic 
$(a_1,a_2,a_3)$ satisfying the condition 
$2 a_1^2 + 2 a_2^2 + a_3^2 \in 5 \Z$;
\begin{equation} \label{twelve1}
\begin{split}
(1,1,1)&, \ (1,1,9), \ (1,9,1), \ (9,1,1), \ (1,3,5), \ (1,7,5),\\
(3,1,5)&, \ (7,1,5), \ (3,3,3), \ (3,3,7),\  (3,7,3), \ (7,3,3), 
\end{split}
\end{equation}
and
\begin{align}\label{twelve2}
\begin{split}
(9,9,9)&, \ (9,9,1), \ (9,1,9), \ (1,9,9), \ (9,7,5), \ (9,3,5), \\
(7,9,5)&, \ (3,9,5), \ (7,7,7), \ (7,7,3), \ (7,3,7), \ (3,7,7),
\end{split}
\end{align}
and $(5,5,5)$. 
\begin{rem}
(1) \ The characteristic $(5,5,5)$ is an odd half integer characteristic 
(see \cite{Mum}), hence $\Theta_{(5,5,5)}(\Omega)$ vanishes identically.
\\ (2) \ By (\ref{chformula}) and (\ref{chformula2}), we see that 
$\Theta_{(a_1,a_2,a_3)}(\Omega)$ is a scalar multiple of 
$\Theta_{(b_1,b_2,b_3)}(\Omega)$ if 
$a_1 + b_1, a_2 + b_2,a_3 + b_3 \in 10 \Z$. So the system in
 (\ref{twelve1}) and the system in (\ref{twelve2}) 
are essentially the same.
\end{rem}
\begin{lem} \label{transitiv}
Let $(a_1,a_2,a_3)$ be a member of the system
 (\ref{twelve1})(equivalently (\ref{twelve2})). The group $\Gamma$ acts
 on the set of twelve $\vartheta(a_1,a_2,a_3)$ transitively.
\end{lem} 
\begin{proof}
We have an explicit form of $\hat{g}_{ij}$ in 
(\ref{matrixg12}) -- (\ref{matrixg45}). We use it and obtain
\[
 \hat{g}_{12}(a_1,a_2,a_3) \equiv (-a_1,a_2,a_3), \quad
\hat{g}_{34}(a_1,a_2,a_3) \equiv (a_1,-a_2,a_3),   
\]
\begin{center}
\begin{tabular}{|c||c|c|}   \hline
$(a_1,a_2,a_3)$ & $\hat{g}_{23}(a_1,a_2,a_3)\equiv$ &
 $\hat{g}_{45}(a_1,a_2,a_3)\equiv$ \\ \hline
(1,1,1)& (3,3,7) & (1,9,9)\\
(1,1,9)& (7,7,7) & (1,7,5) \\
(1,9,1)& (9,7,5) & (1,3,5) \\
(9,1,1)& (7,9,5) & (9,9,9) \\
(1,3,5)& (9,1,9) & (1,9,1) \\
(1,7,5)& (7,3,3) & (1,1,9) \\
(3,1,5)& (1,9,9) & (3,3,7) \\
(7,1,5)& (3,7,3) & (7,3,7) \\
(3,3,3)& (9,9,1) & (3,7,7) \\
(3,3,7)& (1,1,1) & (3,1,5) \\
(3,7,3)& (7,1,5) & (3,9,5) \\
(7,3,3)& (1,7,5) & (7,7,7) \\ \hline
\end{tabular} \end{center}
According to (\ref{formula}),
\[
 g(\vartheta(a_1,a_2,a_3)) = \vartheta(\hat{g}(a_1,a_2,a_3))
\]
So the assertion follows.
\end{proof}
\subsection{The zero loci of twelve theta functions}
Here we state Riemann's theorem. Let $C$ be an algebraic curve of genus 
$g$, let $\{A_i, B_i\}$ be a symplectic basis of $\mathrm{H}_1(C,\Z)$ 
such that $A_i \cdot B_j = \delta_{ij}$, and let $\{\omega_i\}$ be the basis
of $\mathrm{H}^0(C,\Omega^1)$ such that $\int_{A_i} \omega_j =
\delta_{ij}$. Then $\Omega = (\int_{B_i} \omega_j)$ belongs to $\s_g$.
We denotes ${}^t(\int_{\gamma} \omega_1, \cdots, \int_{\gamma}
\omega_g)$ by $\int_{\gamma} \omega$. 
\begin{tm}[see \cite{Mum}, p149] \label{Riemanntm} 
Let us fix a point $P_0 \in C$. Then there is a vector $\Delta \in \C^g$, 
 such that for all $z \in \C^g$, multi-valued function
\[
  f(P) = \Theta_{(0,0)}(z + \int^P_{P_0} \omega,\ \Omega) \quad (P \in C)
\]
on $C$ either vanishes identically, or has $g$ zeros $Q_1, \cdots , Q_g$
 with 
\[
 \sum_{i =1}^g \int_{P_0}^{Q_i} \omega \equiv -z + \Delta 
\quad \text{mod} \  \Omega \Z^g + \Z^g.
\]
\end{tm}
\begin{rem}[see \cite{Mum}] \label{Riemannrem}
(1) \ The vector $\Delta$ in the theorem is called the Riemann constant, 
 and depends on the symplectic basis $\{A_i, b_i\}$ and the base point
 $P_0$. For the fixed $\{A_i, B_i\}$ and $P_0$, $\Delta$ is uniquely 
determined as the point of the Jacobian 
$J(C) = \C^g / (\Omega \Z^g + \Z^g)$ by the property of the theorem.\\
(2) \ If we take $P_0$ such that the divisor $(2g-2)P_0$ is linearly
 equivalent to the canonical divisor, then we have
$\Delta \in \frac{1}{2} \Omega \Z + \frac{1}{2} \Z$.
\end{rem}
\begin{cor}[see \cite{Mum}] \label{Riemanncor}
Under same situation as the theorem, $\Theta_{a,b}(\Omega) = 0$ if and only
 if there exist $Q_1, \cdots, Q_g \in C$ such that
\[
 \Delta - (\Omega a + b) \equiv \sum_{i = 1}^{g-1} \int_{P_0}^{Q_i} \omega.
\]
\end{cor}
\vskip3mm
Now, let us return to our case. Let $\lambda^0 \in \xc$ and $C_0$ be as 
in section \ref{sectionperiod} and $\omega_1, \cdots, \omega_6$ 
be the basis of $\mathrm{H}^0(C_0, \Omega^1)$ 
such that $\int_{A_i} \omega_j = \delta_{ij}$. We denote the ramified 
points over $\lambda_i \in \pr$ by $P_i \in C_0$. Let us take the base
 point $P_0$ arbitrary among $\{P_1, \cdots, P_5\}$ and $\Delta_0$ be the 
Riemann constant with respect to $\{A_i,B_i\}$ and $P_0$.
\begin{lem} \label{lemma(5,5,5)}
The Riemann constant $\Delta_0$ corresponds to the characteristic $(5,5,5)$.
\end{lem}
\begin{proof}
The divisor of the holomorphic 1-form 
$(x-\lambda_i)^2dx/y^4$ is $10 P_i$. Hence $\Delta_0$ is a half integer
 characteristic (see Remark \ref{Riemannrem}). 
For $z = \Omega a + b \ (a,b \in \mathbb{R}^6)$ and
 $\sigma = \begin{pmatrix} A&B \\ C&D \end{pmatrix}$, applying
 (\ref{formula}) we have
\begin{align} \label{equ}
\Theta_{\sigma(\Delta_0 -z)}(\Omega) = 
(\text{a unit function})  \times \Theta_{\Delta_0 -z }(\Omega)
\end{align}
since $\sigma \Omega = \Omega$. By (\ref{gab}) and Remark \ref{Riemanncor}, 
we have 
\[
\sigma(\Delta_0 -z) = \sigma \Delta_0 - {}^t(C \Omega + D)^{-1} z.
\] 
Hence it holds
\begin{align*}
\Theta_{\sigma \Delta_0 - {}^t(C \Omega + D)^{-1}z}(\Omega) = 0 
\ &\Leftrightarrow \
\Theta_{\Delta_0 -z }(\Omega) =0 \\
&\Leftrightarrow \ z \equiv \sum_{i = 1}^{5} \int_{P_0}^{Q_i} \omega 
\quad \text{for} \quad {}^{\exists} Q_1, \cdots, Q_5 \in C_0
\end{align*}
by Corollary \ref{Riemanncor}. Namely, putting 
$w = {}^t(C \Omega + D)^{-1}z$ we have
\begin{align*}
\Theta_{\sigma \Delta_0 - w}(\Omega) = 0 
\ &\Leftrightarrow \
{}^t(C \Omega + D) w \equiv \sum_{i = 1}^{5} \int_{P_0}^{Q_i} \omega 
\quad \text{for} \quad {}^{\exists} Q_1, \cdots, Q_5 \in C_0.
\end{align*}
Let us recall that $\sigma$ is the symplectic representation matrix of
 $\rho$ with respect to the basis $\{A_i, B_i \}$ of 
$\mathrm{H}_1(C_0, \Z)$. And we have
\[
\begin{pmatrix}I & \Omega \end{pmatrix} 
\begin{pmatrix}{}^tD & {}^tB \\ {}^tC & {}^tA\end{pmatrix} = 
\begin{pmatrix}{}^t(C \Omega + D) & {}^t(A \Omega + B)\end{pmatrix} =
{}^t(C \Omega + D) \begin{pmatrix}I & \Omega \end{pmatrix},
\]
so ${}^t(C \Omega + D)$ is the representation matrix of $\rho$ with
 respect to the basis $\{\omega_1, \cdots, \omega_6 \}$ of
 $\mathrm{H}^0(C_0, \Omega^1)$. Hence it holds
\begin{align*}
\Theta_{\sigma \Delta_0 - w}(\Omega) = 0 
\ &\Leftrightarrow \
w \equiv \sum_{i = 1}^{5} \int_{P_0}^{Q_i} (\rho^{-1})^* \omega
\equiv \sum_{i = 1}^{5} \int_{\rho^{-1}(P_0)}^{\rho^{-1}(Q_i)} \omega 
\equiv \sum_{i = 1}^{5} \int_{P_0}^{\rho^{-1}(Q_i)} \omega 
\end{align*}
Recalling Remark \ref{Riemannrem} (1), this implies that 
$\sigma \Delta_0$ is the Riemann constant, that is 
$\sigma \Delta_0 \equiv \Delta_0$. Hence we have $\Delta_0 \equiv (5,5,5)$ 
since $(5,5,5)$ is the unique $\sigma$--invariant half integer
 characteristic.
\end{proof}
\vskip 3mm
Next, let us consider the oriented arcs $\alpha_k(i,j)$ defined by (\ref{arc}) 
and the integrals $\int_{\alpha_k(i,j)} \omega \in \C^6$. 
\vskip3mm
\begin{lem} \label{integrals}
The integral $\int_{\alpha_k(i,j)} \omega$ is a five torsion point
 $\Omega a + b$ on $\C^6 / (\Omega \Z^6 + \Z^6)$ of the form 
\[
 a = \frac{1}{10}{}^t(a_1, a_2, a_3,a_1, a_2, a_3), \
 b = \frac{1}{10}{}^t(-2a_1, -2a_2, -2a_3,-a_1, -a_2, -a_3) 
\]
with $a_1, a_2, a_3 \in 2 \Z$. In explicit way, it holds
\begin{align*}
 \int_{\alpha_k(1,2)} \omega \equiv (6,0,0), \ 
\int_{\alpha_k(1,3)} \omega \equiv (8,2,6), \ 
\int_{\alpha_k(1,4)} \omega \equiv (8,8,6), \ 
\int_{\alpha_k(1,5)} \omega \equiv (8,0,8) \\
\text{mod} \ \Omega \Z^6 + \Z^6
\end{align*}
with the same notation in Definition \ref{defchara} and identification 
referred in Remark \ref{remtrans} (Note that any $\alpha_k(i,j)$ is 
written as a combination of $\alpha_k(1,2)$, $\alpha_k(1,3)$, $\alpha_k(1,4)$ 
 and $\alpha_k(1,5)$).
\end{lem}
\begin{proof}
Since $D_{ij} = \alpha_i(1,5) - \alpha_j(1,5)$ is a cycle, we see that 
$\int_{\alpha_i(1,5)} \omega \equiv \int_{\alpha_j(1,5)} \omega$ mod
 $\Omega \Z^6 + \Z^6$. And we have
\begin{align*}
\int_{D_{12} + D_{15}} \varphi_1 \ =\ 
\int_{2 \alpha_1(1,5) - \alpha_2(1,5) - \alpha_5(1,5)} \varphi_1 \
= \ (2- \z^2 - \z^3) \int_{\alpha_1(1,5)} \varphi_1. 
\end{align*}
By the same calculation, we see that
\begin{align*}
\int_{\alpha_1(1,5)} \varphi_k \ =& \begin{cases}
\frac{1}{5}(2- \z - \z^4) \int_{D_{12} + D_{15}} \varphi_k & (k = 1,2,3) \\
{}\\
\frac{1}{5}(2- \z^2 - \z^3) \int_{D_{12} + D_{15}} \varphi_k & (k = 4,5,6) 
\end{cases} \\
=& \frac{1}{5} \int_{[2 - \rho^2 - \rho^3] (D_{12} + D_{15})} \varphi_k.
\end{align*}
Calculating intersection numbers, we have the following equality
\begin{align*}
 [2 - \rho^2 - \rho^3] (D_{12} + D_{15}) = 
2 A_1 + 2 A_3 + A_4 + A_6 + 4 B_1 + 4 B_3 - B_4 - B_6
\end{align*}
as homology classes. Hence it holds 
\begin{align*}
 \int_{\alpha_1(1,5)} \omega &\equiv \frac{1}{5} 
\int_{2 A_1 + 2 A_3 + A_4 + A_6 + 4 B_1 + 4 B_3 - B_4 - B_6} \omega \\
&\equiv \frac{1}{10}
\int_{-6 A_1 -6 A_3 -8 A_4 -8 A_6 +8 B_1 +8 B_3 +8B_4 +8 B_6} \omega
\equiv (8,0,8).
\end{align*}
By the same way, we obtain the results for $\alpha_k(1,2)$,
$\alpha_k(1,3)$ and $\alpha_k(1,4)$. 
\end{proof} 
\vskip3mm 
Let $C_{\lambda} \ (\lambda \in \xc)$ be any element of our family
$\mathcal{F}$. We defined in Section \ref{sectionperiod} the system
$\{\alpha_k(i,j)(\lambda) \}$, $\{\gamma_i(\lambda) \}$ and 
$\{A_i(\lambda), B_i(\lambda) \}$ on $C_{\lambda}$ depending on the arc
$r$. The point $P_0$ has always the same meaning. So Lemma
\ref{lemma(5,5,5)} and \ref{integrals} are true for $C_{\lambda}$ using
these notations. Let $\Delta \equiv (5,5,5)$ denote the Riemann constant
on $C_{\lambda}$. \\
\indent Now, recall that  $\mathbb{B}_2^0$, $\ell(ij)$ 
stands for $\Phi(\xc)$ and $\pi^{-1}(\tilde{\Phi}(L(ij)))$ 
 respectively(see Section \ref{sectiondl}).
\begin{prop} 
$\vartheta(1,1,1) \cap \mathbb{B}_2^{\circ} = \phi$.
\end{prop}
\begin{proof}
Let us consider a curve $C = C_{\lambda} \ (\lambda \in \xc)$ and its period 
 $\Omega = \Omega_{\lambda}$. We assume that $\Theta_{(1,1,1)}(\Omega) = 0$.
According to Corollary \ref{Riemanncor}, 
there exist points $Q_1, \cdots, Q_5 \in C$ such that 
\[
\sum_{i =1}^5 \int_{P_5}^{Q_i} \omega \equiv \Delta - (1,1,1) 
\equiv (4,4,4). 
\]
On the other hand, by Lemma \ref{integrals}, we have
\[
 \int_{P_4}^{P_3} \omega \equiv (0,4,0), \quad 
\int_{P_5}^{P_1} \omega \equiv (2,0,2).
\]
Hence it holds
\[
 \sum_{i =1}^5 \int_{P_5}^{Q_i} \omega \equiv 
2 \int_{P_5}^{P_1} \omega + \int_{P_4}^{P_3} \omega.
\]
By Abel's theorem, the divisor $\sum_{i =1}^g Q_i$ is linearly
 equivalent to the divisor $D = 2 P_1 + P_3 - P_4 + 3 P_5$, and we have
\begin{align} \label{eq}
 \dim \mathrm{H}^0(C, \mathcal{O}(D)) = 
\dim \mathrm{H}^0(C, \mathcal{O}(\sum_{i =1}^g Q_i)) \geq 1 
\end{align}
For the effective divisor $D' = D + P_0$, we have
\[
\dim \mathrm{H}^0(C, \mathcal{O}(D')) = 
 \dim \mathrm{H}^0(C, \Omega^1(-D')) + 1
\]
by the Riemann-Roch. We claim that 
$\dim \mathrm{H}^0(C, \Omega^1(-D')) = 0$. In fact, the basis 
$\{\varphi_i\}$ is written as
\begin{align*}
 \varphi_1 = y^2 \varphi, \quad \varphi_2 = y \varphi, \quad 
\varphi_3 = x^2 \varphi, \quad \varphi_4 = x \varphi, \quad 
\varphi_5 = xy \varphi, \quad \varphi_5 = \varphi \\
( \varphi = 4 \frac{dy}{f'(x)}, 
\quad f(x) = \prod_{i = 1}^5(x -\lambda_i)),
\end{align*}
and we have following vanishing orders;
\[
 \mathrm{ord}_{P_i}(y) = 1, \quad 
\mathrm{ord}_{P_i}(x - \lambda_j) = 5 \delta_{ij}, \quad
\mathrm{ord}_{P_i}(\varphi) = 0 \quad (i,j = 1, \cdots, 5).
\]
Because any holomorphic 1-form is written in the form
\[
 \text{(inhomogeneous quadratic polynomial of $x,y$)} \times \varphi,
\]
we see that there is no holomorphic 1-form $\xi$ such that
\[
 \mathrm{ord}_{P_1}(\xi) \geq 2, \quad \mathrm{ord}_{P_3}(\xi) \geq 1, \quad
\mathrm{ord}_{P_5}(\xi) \geq 3. 
\]
Hence we have $\dim \mathrm{H}^0(C, \mathcal{O}(D')) = 1$, that is,
$\mathrm{H}^0(C, \mathcal{O}(D'))$ contains only constant functions.
This contradicts to (\ref{eq}) since 
$\mathrm{H}^0(C, \mathcal{O}(D)) \subset \mathrm{H}^0(C, \mathcal{O}(D'))$
and $D$ is not effective.
\end{proof}
\vskip3mm
\begin{cor} \label{cornonvanish}
Let $(a_1,a_2,a_3)$ be a $\sigma$--invariant characteristic in
 (\ref{twelve1}). Then we have
$\vartheta(a_1,a_2,a_3) \cap \mathbb{B}_2^{\circ} = \phi$.
\end{cor}
\begin{proof}
This follows from Lemma \ref{transitiv}.
\end{proof}
\vskip3mm
Hence $\vartheta(a_1,a_2,a_3)$ is the union of certain $\ell(ij)$'s.
\begin{lem} \label{evaluation}
Let $\eta_0$ be the point $[0:0:1] \in \B$, and let $(a_1,a_2,a_3)$ be a 
 member of (\ref{twelve1}). If $a_1, a_2,a_3 \in \{1,9\}$, then we have 
$\Theta_{(a_1,a_2,a_3)}(\Omega(\eta_0)) \ne 0$. 
\end{lem}
\begin{proof}
Let 
\[
(a',b') = (\frac{1}{10}{}^t(\alpha,\alpha),
\frac{1}{10}{}^t(-2 \alpha,-\alpha))
\]
be a characteristic in $(\mathbb{Q}^2)^2$. Let $\Theta_{\alpha}(\tau)$
 denote the theta constant $\Theta_{(a',b')}(\tau) \ (\tau \in \s_2)$.
Using this notation, we have
\begin{align} \label{3theta}
 \Theta_{(a_1,a_2,a_3)}(\Omega(\eta_0)) = 
\Theta_{a_1}(\tau_0) \Theta_{a_2}(\tau_0) \Theta_{a_3}(\tau_0), \quad
\tau_0 = \begin{pmatrix} \z - 1 &\z^2 +\z^3 \\ 
\z^2 +\z^3 & -\z^4 \end{pmatrix}
\end{align}
(see (\ref{eta_0})). So our assertion is reduced to the inequality 
$\Theta_1(\tau_0) \ne 0$, since $\Theta_9$ is a constant multiple of 
$\Theta_1$. Set
\[
 a = {}^t(\frac{1}{10},\frac{1}{10}), \quad 
b = {}^t(-\frac{2}{10},-\frac{1}{10}), \quad n = {}^t(n_1,n_2),
\]
and set
\[
 f(n_1,n_2) = \exp[\pi \sqrt{-1} (^t(n+a) \tau_0 (n+a) + 2^t(n+a)b)].
\]
By definition, $\Theta_1(\tau_0) = \sum_{n_1,n_2 \in \Z} f(n_1,n_2)$.
For simplicity, we denote $n + a$ by $m = (m_1,m_2)$. 
By elementary calculations, we see that
\[
 |f(n_1,n_2)| = 
\exp[-\pi \sin(\frac{2 \pi}{5})\{m_1^2  + (3 - \sqrt{5})m_1 m_2 + m_2^2\}].
\]
In case $m_1 m_2 >0$, we have
\[
 |f(n_1,n_2)| < \exp[-\pi \sin(\frac{2 \pi}{5})\{m_1^2  + m_2^2\}].
\]
In case $m_1 m_2 < 0$, we have
\begin{align*}
 |f(n_1,n_2)| &< \exp[-\pi \sin(\frac{2 \pi}{5})\{m_1^2 +m_1 m_2
 +m_2^2\}] \\
&= \exp[-\pi \sin(\frac{2 \pi}{5})\{\frac{1}{2}(m_1^2 + m_2^2) + 
\frac{1}{2}(m_1 + m_2)^2\}] \\
&< \exp[-\frac{\pi}{2} \sin(\frac{2 \pi}{5})\{m_1^2 + m_2^2\}].
\end{align*}
Consequently, 
\[
 |f(n_1,n_2)| < \alpha^{m_1^2 + m_2^2}, \quad 
(\alpha = \exp[-\frac{\pi}{2} \sin(\frac{2 \pi}{5})])
\]
for any $n_1, n_2 \in \Z$. Set
\[
 D_1 = \{ (n_1,n_2) \in \Z^2 \ : \ -10 \leq n_1, n_2 \leq 10\}, \quad
D_2 = \Z^2 - D_1,
\]
and consider the summations
\[
 S_1 = \sum_{D_1} f(n_1,n_2), \quad S_2 = \sum_{D_2} f(n_1,n_2).
\]
Using a computer, we can evaluate $|S_1|$ and $|S_2|$.
We have a approximate value
\[
 |S_1| \fallingdotseq 1.13746 \cdots,
\]
by {\it Mathematica}. On the other hand, we have
\[
|S_2| < \sum_{D_2}|f(n_1,n_2)| < \sum_{D_2} \alpha^{m_1^2 + m_2^2}.
\]
The last term is very small. For example,
\[
 \sum_{n_1 \geq 10, n_2 \geq 0} \alpha^{m_1^2 + m_2^2} <
(\sum_{n_1 \geq 10} \alpha^{n_1})(\sum_{n_2 \geq 0} \alpha^{n_2}) =
(\frac{\alpha^{10}}{1-\alpha})( \frac{\alpha}{1-\alpha})
 \fallingdotseq 5.40545 \times 10^{-7},
\]
and the same calculations shows $|S_1| \gg |S_2|$. 
This implies $\Theta_1(\tau_0) = S_1 + S_2 \ne 0$.
\end{proof}
\vskip3mm
\begin{lem} \label{vanishinglemma}
(1) \ If we have $a_1 \equiv 3,7 \ \mathrm{mod}\ 10$, then 
$\Theta_{(a_1,a_2,a_3)}$ vanishes on $\ell(12)$. \\
(2) \ If we have $a_2 \equiv 3,7 \ \mathrm{mod}\ 10$, then
$\Theta_{(a_1,a_2,a_3)}$ vanishes on $\ell(34)$.
\end{lem}
\begin{proof}
Set $g = \hat{g}_{12} = \begin{pmatrix} A&B\\C&D\end{pmatrix}$ and set
$\Omega = \Omega(\eta)$ with $\eta = [0:\eta_2:\eta_3] \in \B$.
By the computation same as the one in the proof of Lemma
 \ref{vanish100}, we have
\[
 g \Omega = \Omega, \quad \det(C \Omega + D) = \z, \quad
\phi_{(a_1,a_2,a_3)}(g) = \frac{1}{40} a_1^2, \quad
\Theta_{g(a_1,a_2,a_3)}(\Omega) = \Theta_{(a_1,a_2,a_3)}(\Omega).
\]
Hence it holds
\[
\Theta_{(a_1,a_2,a_3)}(\Omega)^8 = 
\exp[\frac{2}{5} \pi \sqrt{-1}(a_1^2 - 1)] 
\Theta_{(a_1,a_2,a_3)}(\Omega)^8 
\]
Therefore $\Theta_{(a_1,a_2,a_3)}(\Omega(\eta))$ vanishes on 
the mirror of $g_{12}$ provided $a_1 \equiv 3,7 \ \mathrm{mod}\ 10$. 
This implies assertion (1). The assertion (2) follows by the same
 argument with $g = g_{34}$ and $\eta = [\eta_1: 0:\eta_3] \in \B$.
\end{proof}
\begin{prop} \label{proptable}
We have the Table \ref{table} for the vanishing loci of 
twelve theta constants coming from the system (\ref{twelve1}). 
In the table, ``v'' implies that $\Theta_{(a_1,a_2,a_3)}$ vanishes
 there, and the blank implies $\Theta_{(a_1,a_2,a_3)}$ is not
 identically zero there. For example, $\Theta_{(1,1,1)}$ vanishes on 
$\ell(13)$ and is not identically zero on $\ell(12)$.
\begin{table}[hbtp]
\begin{center}
\begin{tabular}{|c||c|c|c|c|c|c|c|c|c|c|}   \hline
$(a_1,a_2,a_3)$ & $\ell(12)$ & $\ell(13)$ & $\ell(14)$ & 
$\ell(15)$ & $\ell(23)$ & $\ell(24)$ & $\ell(25)$ & $\ell(34)$ & 
$\ell(35)$ & $\ell(45)$ \\ \hline
$(1,1,1)$ &  &v &  &v &v &v &  &  &  &v \\ \hline
$(1,1,9)$ &  &v &v &  &  &v &v &  &v &  \\ \hline
$(1,9,1)$ &  &  &v &v &v &v &  &  &v &  \\ \hline
$(9,1,1)$ &  &v &v &  &v &  &v &  &  &v \\ \hline
$(1,3,5)$ &  &  &v &v &v &  &v &v &  &  \\ \hline
$(1,7,5)$ &  &v &  &v &  &v &v &v &  &  \\ \hline
$(3,1,5)$ &v &  &v &  &v &  &  &  &v &v \\ \hline
$(7,1,5)$ &v &v &  &  &  &v &  &  &v &v \\ \hline
$(3,3,3)$ &v &  &v &  &  &  &v &v &v &  \\ \hline
$(3,3,7)$ &v &  &  &v &v &  &  &v &  &v \\ \hline
$(3,7,3)$ &v &v &  &  &  &  &v &v &  &v \\ \hline
$(7,3,3)$ &v &  &  &v &  &v &  &v &v &  \\ \hline
\end{tabular} \end{center} \caption{} \label{table}
\end{table}
\end{prop}
\begin{proof}
By Lemma \ref{vanishinglemma}, 
\[
 \Theta_{(3,1,5)}, \quad \Theta_{(7,1,5)}, \quad \Theta_{(3,3,3)}, \quad 
\Theta_{(3,3,7)}, \quad \Theta_{(3,7,3)}, \quad \Theta_{(7,3,3)}  
\]
vanish on $\ell(12)$, and
\[
 \Theta_{(1,3,5)}, \quad \Theta_{(1,3,5)}, \quad \Theta_{(3,3,3)}, \quad 
\Theta_{(3,3,7)}, \quad \Theta_{(3,7,3)}, \quad \Theta_{(7,3,3)}  
\]
vanish on $\ell(34)$. By Lemma \ref{evaluation},
\[
 \Theta_{(1,1,1)}, \quad \Theta_{(1,1,9)}, \quad \Theta_{(1,9,1)}, \quad 
\Theta_{(9,1,1)}
\]
are not identically zero on $\ell(12)$ and on $\ell(34)$, since 
$\eta_0 = [0:0:1] \in \ell(12) \cap \ell(34)$. The result is obtained by 
 applying the transformation formula (\ref{formula}) for above theta
 constants and $\hat{g}_{ij}$. For example, we have
\[
 \Theta_{\hat{g}_{12}\hat{g}_{45}(a_1,a_2,a_3)}
(\hat{g}_{12}\hat{g}_{45} \Omega) =
\text{(a unit function)} \times \Theta_{(a_1,a_2,a_3)}(\Omega).
\]
Since $\hat{g}_{12}\hat{g}_{45}(1,3,5) \equiv (9,9,1)$ 
(see Lemma \ref{transitiv}) and
$g_{12}g_{45}(\ell(12)) = \ell(12)$, we see that $\Theta_{(1,3,5)}$ is
not identically zero on $\ell(12)$.           
\end{proof}
\subsection{Automorphic Factor}
We study the automorphic factor appeared in the transformation
formula (\ref{formula}) with respect to $\Gamma(1-\z)$ and $\Omega =
\Omega(\eta) $. Let $H$ be the diagonal matrix 
$\mathrm{diag}(1,1,-\z^3(1 + \z))$. We denote ${}^t\eta H \eta$ by
$\left<\eta,\eta\right>$. Set
\[
 F_g(\eta) = \frac{\left<g\eta,g\eta\right>}{\left<\eta,\eta\right>}
\]
for $g \in \Gamma$ and $\eta \in \B$. Obviously, we have 
the following lemma.
\begin{lem} \label{lemcocycle}
$F_g(\eta)$ satisfies the cocycle condition with respect to $\Gamma$. That
 is,
\[
 F_{g_1g_2}(\eta) = F_{g_1}(g_2\eta) F_{g_2}(\eta),
\qquad g_1, g_2 \in \Gamma.
\]
\end{lem}
\vskip3mm
\begin{prop}
There exist the non trivial character 
\[
 \chi \ :\  \Gamma \longrightarrow \mu_5 = \{1,\z,\cdots, \z^4\}
\]
such that 
\begin{align*}
 \det(C \Omega(\eta) + D) = \chi(g) F_g(\eta) \quad (\eta \in \B)
\end{align*}
for $g \in \Gamma$, where the matrix
$\begin{pmatrix} A&B\\C&D\end{pmatrix}$ is the symplectic representation 
 $\hat{g}$ of $g$. 
\end{prop}
\begin{proof}
According to the case by case calculation, we have
\[
\det(C \Omega(\eta) + D) = \z^3 F_g(\eta) \quad (\eta \in \B)
\]
for $g = g_{12}, g_{23}, g_{34}, g_{45}$. Since 
$\det(C \Omega(\eta) + D)/F_g(\eta)$ satisfies the cocycle condition, we
 obtain the result.
\end{proof}
Now let $(a,b)$ be a invariant characteristic $(a_1,a_2,a_3)$, and
$(a_g,b_g)$ be $\hat{g}(a,b)$ for $g \in \Gamma(1-\z)$. Since
\[
(a_g,b_g) \equiv (a,b) \quad \text{mod} \ \Z,
\] 
we have
\[
 \Theta_{\hat{g}(a,b)}(\Omega) = \Theta_{(a_g,b_g)}(\Omega) 
= \Theta_{(a_g - a + a,b_g - b + b)}(\Omega)
= \exp[2 \pi \sqrt{-1} {}^ta(b_g - b)] \Theta_{(a,b)}(\Omega)
\]
by (\ref{chformula}). Set
\[
 \phi'_{(a_1,a_2,a_3)}(\hat{g}) = 
\phi_{(a_1,a_2,a_3)}(\hat{g}) - {}^ta(b_g - b).
\]
Then we can write the transformation formula (\ref{formula}) as
\begin{align} \label{newformula}
\Theta_{(a_1,a_2,a_3)}(\Omega(g\eta)) = 
\kappa(\hat{g})\exp(2 \pi \sqrt{-1}\phi'_{(a_1,a_2,a_3)}(\hat{g}))
[\chi(g) F_g(\eta)]^{\frac{1}{2}}\Theta_{(a_1,a_2,a_3)}(\Omega(\eta)),
\end{align}
where $\kappa(\hat{g})$ is a 8-th root of 1 depending only on $\hat{g}$.
\vskip3mm
\begin{lem} \label{lemmaphi'}
Let $g$ be in $\Gamma(1 - \z)$. Then, the values
\[ 
[\exp(2 \pi \sqrt{-1}\phi'_{(a_1,a_2,a_3)}(\hat{g}))]^5 
\]
are the same for all twelve characteristics $(a_1,a_2,a_3)$ 
in (\ref{twelve1}).
\end{lem}
\begin{proof}
By direct calculations, we have
\begin{align*}
5 \phi'_{(a_1,a_2,a_3)}(\hat{h}_{12}) \equiv \frac{1}{8}, \quad
5 \phi'_{(a_1,a_2,a_3)}(\hat{h}_{13}) \equiv \frac{3}{4}, \quad
5 \phi'_{(a_1,a_2,a_3)}(\hat{h}_{14}) \equiv \frac{1}{2}, \\
5 \phi'_{(a_1,a_2,a_3)}(\hat{h}_{23}) \equiv \frac{1}{2}, \quad
5 \phi'_{(a_1,a_2,a_3)}(\hat{h}_{34}) \equiv \frac{3}{4} \quad
(\text{mod} \ \Z)
\end{align*}
for the twelve $(a_1,a_2,a_3)$. According to Lemma \ref{lemcocycle}, the 
 equality (\ref{newformula}) shows that 
\[
\kappa(\hat{g})\exp[2 \pi \sqrt{-1}\phi'_{(a_1,a_2,a_3)}(\hat{g})] 
\]
is a character on $\Gamma(1 -\z)$. So we obtain the
 result for any $g \in \Gamma(1-\z)$.
\end{proof}
\vskip3mm
\begin{cor}
Let $(a_1, a_2, a_3)$ and $(b_1, b_2, b_3)$ be in (\ref{twelve1}). Then, 
 the function 
\[
 \frac{\Theta_{(a_1, a_2, a_3)}(\Omega(\eta))^5}
{\Theta_{(b_1, b_2, b_3)}(\Omega(\eta))^5}
\]
is well-defined as meromorphic function on $\B / \Gamma(1-\z)$.
\end{cor}
\vskip3mm
Let $\Omega = \Omega_{\lambda}$ be the period matrix of a curve
$C_{\lambda} \ (\lambda \in \xc)$, $P_0$  be a ramified point of $C
\rightarrow \pr$.
\vskip 3mm
\begin{prop} 
Let $(a_1, a_2, a_3)$ and $(b_1, b_2, b_3)$ be in (\ref{twelve1}). 
The function 
\[
 f(P) = \frac{\Theta_{(a_1, a_2, a_3)}
(\int_{P_0}^{P} \omega, \ \Omega)^5}
{\Theta_{(b_1, b_2, b_3)}(\int_{P_0}^{P} \omega, \ \Omega)^5}
\quad (P \in C_{\lambda})
\]
is a single-valued meromorphic function on $C_{\lambda}$, where the
 paths of integrations in the numerator and the denominator are 
chosen as  same.      
\end{prop}
\begin{proof}
Note that Corollary \ref{cornonvanish} asserts 
\[
 \Theta_{(a_1, a_2, a_3)}(\int_{P_0}^{P_0} \omega, \ \Omega)
= \text{const.} \times \Theta_{(a_1, a_2, a_3)}(0, \ \Omega)
\ne 0,
\]
where the constant depends on the path of integration.
So the numerator is not identically zero, and it is same for the
 denominator. By the assumption we have
\[
 (a_1, a_2, a_3) - (b_1, b_2, b_3) \in (\frac{1}{5} \Z^6)^2.
\]
By using the formula (\ref{translate1}) and (\ref{translate2}) 
we can check that
\[
 \frac{\Theta_{(a_1, a_2, a_3)}
(\int_{P_0}^{P} \omega + \Omega m + n, \ \Omega)^5}
{\Theta_{(b_1, b_2, b_3)}
(\int_{P_0}^{P} \omega + \Omega m + n, \ \Omega)^5} =
\frac{\Theta_{(a_1, a_2, a_3)}(\int_{P_0}^{P} \omega, \ \Omega)^5}
{\Theta_{(b_1, b_2, b_3)}(\int_{P_0}^{P} \omega, \ \Omega)^5}
\]
for $m, n \in \Z^6$. This implies single-valuedness of $f$.
\end{proof}
\vskip3mm
Let us consider the meromorphic function  
\[
f(P) = \frac{\Theta_{(1,1,1)}(\int_{P_1}^{P} \omega, \ \Omega)^5}
{\Theta_{(3,3,7)}(\int_{P_1}^{P} \omega, \ \Omega)^5}
\]
on $C_{\lambda}$. By Lemma \ref{integrals}, we have
\begin{align*}
\Delta - (1,1,1) \equiv (4,4,4) \equiv 2 \int_{P_1}^{P_2} \omega + 
3 \int_{P_1}^{P_3} \omega + \int_{P_1}^{P_4} \omega, \\
\Delta - (3,3,7) \equiv (2,2,8) \equiv 3 \int_{P_1}^{P_2} \omega + 
2 \int_{P_1}^{P_3} \omega + \int_{P_1}^{P_4} \omega.
\end{align*}
By Corollary \ref{Riemanncor},  the zero divisor of 
$\Theta_{(1,1,1)}(\int_{P_1}^{P} \omega, \ \Omega)$ and 
$\Theta_{(3,3,7)}(\int_{P_1}^{P} \omega, \ \Omega)$
are $2 P_2 + 3 P_3 + P_4$ and $3 P_2 + 2 P_3 + P_4$ respectly. Hence we
can write 
\[
 f(P) = c \ \frac{x(P) - \lambda_3}{x(P) - \lambda_2},
\]
where $x(P)$ is the coordinate function 
$x \in \C[x,y]/(y^5 - \prod(x-\lambda_i))$ and $c \ne 0$ is a certain
constant. By Lemma \ref{integrals}, 
\[
 \int_{P_1}^{P_1} \omega \equiv (0,0,0), \quad
\int_{P_1}^{P_5} \omega \equiv (8,0,8).
\]
Substitutes $P = P_1, P_5$ in the above form, then we obtain
\[
\frac{\Theta_{(1,1,1)}((0,0,0), \ \Omega)^5}
{\Theta_{(3,3,7)}((0,0,0), \ \Omega)^5} = 
c \ \frac{\lambda_1 - \lambda_3}{\lambda_1 - \lambda_2}, \quad
\frac{\Theta_{(1,1,1)}((8,0,8), \ \Omega)^5}
{\Theta_{(3,3,7)}((8,0,8), \ \Omega)^5} = 
c \ \frac{\lambda_5 - \lambda_3}{\lambda_5 - \lambda_2}.
\]
Set $(8,0,8) = \Omega \varepsilon' + \varepsilon''$. 
By elementary and patient calculation, we have
\begin{align*}
\Theta_{(1,1,1)}((8,0,8), \ \Omega)^5 = - \z^2
\exp[-5 \pi \sqrt{-1}{}^t\varepsilon' \Omega \varepsilon' -10 \pi \sqrt{-1}
{}^t\varepsilon' \varepsilon''] \Theta_{(1,9,1)}(\Omega)^5 \\
\Theta_{(3,3,7)}((8,0,8), \ \Omega)^5 = 
\exp[-5 \pi \sqrt{-1}{}^t\varepsilon' \Omega \varepsilon' -10 \pi \sqrt{-1}
{}^t\varepsilon' \varepsilon''] \Theta_{(1,3,5)}(\Omega)^5 
\end{align*}
Eliminating $c$, we have the following equality
\[
 \frac{\Theta_{(1,1,1)}(\Omega)^5 \Theta_{(1,3,5)}(\Omega)^5}
{\Theta_{(3,3,7)}(\Omega)^5 \Theta_{(1,9,1)}(\Omega)^5} = 
-\z^2 \frac{(\lambda_1 - \lambda_3)(\lambda_5 - \lambda_2)}
{(\lambda_1 - \lambda_2)(\lambda_5 - \lambda_3)}.
\]
Note that we can regard the above equality as that of meromorphic
functions on $\B / \Gamma(1-\z) \cong \x$.
By the above equality and Proposition \ref{proptable}, we see that 
\begin{enumerate}
\item The vanishing order of $\Theta_{(1,1,1)}(\Omega(\eta))^5$ on
      $\tilde{\Phi}(L(13))$ is 1,
\item The vanishing order of $\Theta_{(1,3,5)}(\Omega(\eta))^5$ on
      $\tilde{\Phi}(L(25))$ is 1,
\item The vanishing order of $\Theta_{(3,3,7)}(\Omega(\eta))^5$ on
      $\tilde{\Phi}(L(12))$ is 1,
\item The vanishing order of $\Theta_{(1,9,1)}(\Omega(\eta))^5$ on
      $\tilde{\Phi}(L(35))$ is 1.
\end{enumerate}
Because $\Gamma$ acts transitively on the set of
$\sigma$--invariant characteristics (see Lemma \ref{transitiv}), we obtain
the following result.
\vskip3mm
\begin{prop} \label{proporder}
Let $(a_1,a_2,a_3)$ be a $\sigma$--invariant characteristic.
If the multi-valued function $\Theta_{(a_1,a_2,a_3)}(\Omega(\eta))^5$ on 
$\B / \Gamma(1-\z)$ identically vanishes on 
$\tilde{\Phi}(L(ij)) = \ell(ij) / \Gamma(1-\z)$, 
then the vanishing order is 1. 
\end{prop}
\vskip7mm
\section{Conclusion}
Now we state our results. \par 
\noindent
{\bf ----- The Schwarz inverse for the Appell HGDE
$F_1(\frac{3}{5},\frac{3}{5},\frac{2}{5},\frac{6}{5})$ ----}
\par
Recall the embedding of $J: \x \rightarrow \mathbb{P}^{11}$ with 
\[
 J(ijklm)(\lambda) = (\lambda_i - \lambda_j)(\lambda_j - \lambda_k)
(\lambda_k - \lambda_l)(\lambda_l - \lambda_m)(\lambda_m - \lambda_i)
\]
in Proposition \ref{propJ}, and the extended period map $\tilde{\Phi}$ 
in section \ref{sectiondl}. 
\begin{tm} \label{maintheorem}
We have a commutative diagram:

\begin{figure}[htbp] \begin{center}
\setlength{\unitlength}{1mm}
\begin{picture}(60,30)
\put(13,28){\vector(1,0){40}}
\put(7,25){\vector(0,-1){20}} 
\put(50,25){\vector(-2,-1){40}} 
\put(0,27){$\x$} \put(30,30){$\tilde{\Phi}$} 
\put(54,27){$\B / \Gamma(1-\z)$} 
\put(5,0){$\mathbb{P}^{11}$} \put(3,16){$J$} \put(30,10){$\varTheta$}
\end{picture} 
\end{center} 
\caption{} \label{diagram}
\end{figure}

by putting
\begin{align} \label{exactmap}
\varTheta = \begin{bmatrix}
\Theta_{(1,1,1)}(\Omega(\eta))^5 \\ \Theta_{(1,1,9)}(\Omega(\eta))^5 \\
\Theta_{(1,9,1)}(\Omega(\eta))^5 \\ \Theta_{(9,1,1)}(\Omega(\eta))^5 \\
\Theta_{(1,3,5)}(\Omega(\eta))^5 \\ \Theta_{(1,7,5)}(\Omega(\eta))^5 \\
\Theta_{(3,3,3)}(\Omega(\eta))^5 \\ \Theta_{(3,3,7)}(\Omega(\eta))^5 \\
\Theta_{(3,7,3)}(\Omega(\eta))^5 \\ \Theta_{(7,3,3)}(\Omega(\eta))^5 \\
\Theta_{(7,1,5)}(\Omega(\eta))^5 \\ \Theta_{(3,1,5)}(\Omega(\eta))^5 
\end{bmatrix}, \quad
J = \begin{bmatrix}
c_1 J(13245)(\lambda) \\ c_2 J(13524)(\lambda) \\ c_3 J(15324)(\lambda) 
\\ c_4 J(13254)(\lambda) \\ c_5 J(15234)(\lambda) \\ c_6 J(13425)(\lambda)
\\ d_1 J(12534)(\lambda) \\ d_2 J(12345)(\lambda) \\ d_3 J(13452)(\lambda)
\\ d_4 J(15342)(\lambda) \\ d_5 J(12453)(\lambda) \\ d_6 J(12354)(\lambda)
\end{bmatrix},
\end{align}
with constants
\[
 [c_1 : \cdots : c_6 : d_1 : \cdots : d_6] =
[1:-1:1:1:\z^3:\z^3:-\z:\z:\z:-\z:-1:-1] \in \mathbb{P}^{11}.
\]
Moreover the map $\varTheta$ is an embedding.
\end{tm}
\begin{proof}
By Proposition \ref{proptable} and Proposition \ref{proporder},
the zero divisor of the i-th component of $\varTheta$ coincides with that of
 the i-th component of $J$ via the isomorphism $\tilde{\Phi}$. 
So we can write as 
\begin{align}
\frac{\Theta_{(1,1,1)}(\Omega(\tilde{\Phi}(\lambda)))^5}
{J(13245)(\lambda)} = c_1 ,\ \cdots \cdots, \
\frac{\Theta_{(3,1,5)}(\Omega(\tilde{\Phi}(\lambda)))^5}
{J(12354)(\lambda)} = d_6 
\end{align}
with certain constants $c_1, \cdots, d_6$. It shows the diagram in
 question is commutative. 
Since $J$ is an embedding and $\tilde{\Phi}$ is an isomorphism, 
we see that $\varTheta$ is an embedding.
\par
\indent Next we determine the ratios of the constants $c_i, d_i$.
\par
\noindent Let $\eta \in \B$ be a point of the form $[0:\eta_2:\eta_3]$. 
For such $\eta$, we have the decomposition
\[
\Omega(\eta) = \tau_0 \oplus \Omega(\eta)'
\] 
(see (\ref{degeperiod})). So we have the splitting
\[
\Theta_{(a_1,a_2,a_3)}(\Omega(\eta)) = \Theta_{a_1} (\tau_0) 
\Psi_{(a_2,a_3)}(\Omega(\eta)'),
\]
where $\Theta_{a_1} (\tau_0)$ is same as in the proof of Lemma 
\ref{evaluation}, and $\Psi_{(a_2,a_3)}(\Omega(\eta)')$ is the theta
 function $\Theta_{(a,b)}(\Omega(\eta)')$ of degree 4 with the characteristic
\[
 a = \frac{1}{10}(a_2,a_3,a_2,a_3), \quad 
b = \frac{1}{10}(-2 a_2,-2 a_3,- a_2,- a_3).
\]
By the same way, we have 
\[
\Theta_{(a_1,a_2,a_3)}(\Omega(\eta)) = \Theta_{a_2} (\tau_0) 
\Psi_{(a_1,a_3)}(\Omega(\eta)'') 
\]
for $\eta = [\eta_1 : 0 : \eta_3] \in \B$. We can check that
\[
\Theta_{1} (\tau_0)^5 = - \Theta_{9} (\tau_0)^5, \quad
\Psi_{(1,9)}(\Omega(\eta)')^5 = \Psi_{(9,1)}(\Omega(\eta)')^5, \quad
\Psi_{(3,5)}(\Omega(\eta)')^5 = \Psi_{(7,5)}(\Omega(\eta)')^5
\] 
by (\ref{chformula}) and (\ref{chformula2}). So we have
\begin{align} \label{ident1}
\frac{\Theta_{(1,1,1)}(\Omega(\eta))^5}{\Theta_{(9,1,1)}(\Omega(\eta))^5}
= \frac{\Theta_1(\tau_0)^5 \Psi_{(1,1)}(\Omega(\eta)')^5}
{\Theta_9(\tau_0)^5 \Psi_{(1,1)}(\Omega(\eta)')^5} = -1
\quad \text{on} \ \ell(12). 
\end{align}
By the same way, we see
\begin{align} \label{ident2}
\frac{\Theta_{(1,1,9)}(\Omega(\eta))^5}{\Theta_{(1,9,1)}(\Omega(\eta))^5}
= 1, \quad
\frac{\Theta_{(1,3,5)}(\Omega(\eta))^5}{\Theta_{(1,7,5)}(\Omega(\eta))^5}
= 1 \quad \text{on} \ \ell(12),
\end{align}
and 
\begin{align} \label{ident3}
\frac{\Theta_{(1,1,1)}(\Omega(\eta))^5}{\Theta_{(1,9,1)}(\Omega(\eta))^5}
= -1, \quad
\frac{\Theta_{(9,1,1)}(\Omega(\eta))^5}{\Theta_{(1,1,9)}(\Omega(\eta))^5}
= 1, \quad
\frac{\Theta_{(3,1,5)}(\Omega(\eta))^5}{\Theta_{(7,1,5)}(\Omega(\eta))^5}
= 1 \quad \text{on} \ \ell(34).
\end{align}
Moreover, we have
\begin{align} \label{ident4}
\frac{\Theta_{(1,1,9)}(\Omega(\eta))^5}{\Theta_{(1,1,1)}(\Omega(\eta))^5} 
= \frac{\Theta_1(\tau_0)^5 \Theta_1(\tau_0)^5 \Theta_9(\tau_0)^5 }
{\Theta_1(\tau_0)^5 \Theta_1(\tau_0)^5 \Theta_1(\tau_0)^5}
= -1  \quad for \ \eta \in \ell(12) \cap \ell(34)
\end{align}
(see (\ref{3theta})). By the transformation formula (\ref{formula}),
 we have
\begin{align*}
\frac{\Theta_{\hat{g}(a_1,a_2,a_3)}(\Omega(g \eta))^5}
{\Theta_{\hat{g}(b_1,b_2,b_3)}(\Omega(g \eta))^5}
= \exp[2 \pi \sqrt{-1}\{\phi_{(a_1,a_2,a_3)}(\hat{g}) - 
\phi_{(b_1,b_2,b_3)}(\hat{g})\}] \frac{\Theta_{(a_1,a_2,a_3)}(\Omega(\eta))^5}
{\Theta_{(b_1,b_2,b_3)}(\Omega(\eta))^5}
\end{align*}
for any pair of $\sigma$--invariant characteristics 
$(a_1,a_2,a_3),(b_1,b_2,b_3)$, and $g \in \Gamma$. By explicit
 calculation of the above formula, we obtain
\begin{align*}
\frac{\Theta_{(1,1,1)}(\Omega(\eta))^5}{\Theta_{(9,1,1)}(\Omega(\eta))^5}
&=& - \z^4 \frac{\Theta_{(3,3,7)}(\Omega(g_{23} \eta))^5}
{\Theta_{(3,1,5)}(\Omega(g_{23} \eta))^5}, \quad
\frac{\Theta_{(1,1,9)}(\Omega(\eta))^5}{\Theta_{(1,9,1)}(\Omega(\eta))^5}
&=& \z^2 \frac{\Theta_{(3,3,3)}(\Omega(g_{23} \eta))^5}
{\Theta_{(1,3,5)}(\Omega(g_{23} \eta))^5}, \\
\frac{\Theta_{(1,3,5)}(\Omega(\eta))^5}{\Theta_{(1,7,5)}(\Omega(\eta))^5}
&=& \z \frac{\Theta_{(1,9,1)}(\Omega(g_{23} \eta))^5}
{\Theta_{(7,3,3)}(\Omega(g_{23} \eta))^5}, \quad
\frac{\Theta_{(3,1,5)}(\Omega(\eta))^5}{\Theta_{(7,1,5)}(\Omega(\eta))^5}
&=& \z \frac{\Theta_{(9,1,1)}(\Omega(g_{23} \eta))^5}
{\Theta_{(3,7,3)}(\Omega(g_{23} \eta))^5}, \\
\frac{\Theta_{(3,1,5)}(\Omega(\eta))^5}{\Theta_{(7,1,5)}(\Omega(\eta))^5}
&=& - \frac{\Theta_{(3,3,7)}(\Omega(g_{45} \eta))^5}
{\Theta_{(3,7,3)}(\Omega(g_{45} \eta))^5}, \quad
\frac{\Theta_{(1,1,9)}(\Omega(\eta))^5}{\Theta_{(1,1,1)}(\Omega(\eta))^5}
&=& \z^2 \frac{\Theta_{(1,7,5)}(\Omega(g_{45} \eta))^5}
{\Theta_{(9,1,1)}(\Omega(g_{45} \eta))^5}.
\end{align*}
Comparing these with  (\ref{ident1}), (\ref{ident2}), (\ref{ident3}) 
and (\ref{ident4}), we have
\begin{align} \label{ident5}
\frac{\Theta_{(3,3,7)}(\Omega(\eta))^5}{\Theta_{(3,1,5)}(\Omega(\eta))^5}
= \z, \quad
\frac{\Theta_{(3,3,3)}(\Omega(\eta))^5}{\Theta_{(1,3,5)}(\Omega(\eta))^5}
= \z^3 \quad 
\frac{\Theta_{(1,9,1)}(\Omega(\eta))^5}{\Theta_{(7,3,3)}(\Omega(\eta))^5}
= \z^4 \quad \text{on} \ \ell(13),
\end{align}
\begin{align} \label{ident6}
\frac{\Theta_{(9,1,1)}(\Omega(\eta))^5}{\Theta_{(3,7,3)}(\Omega(\eta))^5}
= \z^4 \quad \text{on} \ \ell(24),
\end{align}
\begin{align} \label{ident7}
\frac{\Theta_{(3,3,7)}(\Omega(\eta))^5}{\Theta_{(3,7,3)}(\Omega(\eta))^5}
= -1, \quad \text{on} \ \ell(35),
\end{align}
\begin{align} \label{ident8}
\frac{\Theta_{(1,7,5)}(\Omega(\eta))^5}{\Theta_{(9,1,1)}(\Omega(\eta))^5}
= \z^3 \quad \text{on} \ \ell(12) \cap \ell(35),
\end{align}
since it holds
\begin{align*}
g_{23}(\ell(12)) = \ell(13), \quad g_{23}(\ell(34)) = \ell(24), \quad 
g_{45}(\ell(12)) = \ell(12), \quad g_{45}(\ell(34)) = \ell(35).
\end{align*}
Because the commutativity of diagram is established, by using (\ref{ident1}),
 we see that
\begin{align*}
-1 = \frac{\Theta_{(1,1,1)}(\Omega(\eta))^5}{\Theta_{(9,1,1)}(\Omega(\eta))^5}
|_{\ell(12)} &= \frac{c_1}{c_4} \frac{J(13245)}{J(13254)}|_{L(12)} \\
&= \frac{c_1}{c_4} \frac{(\lambda_1 - \lambda_3)(\lambda_3 - \lambda_1)
(\lambda_1 - \lambda_4)(\lambda_4 - \lambda_5)(\lambda_5 - \lambda_1)}
{(\lambda_1 - \lambda_3)(\lambda_3 - \lambda_1)(\lambda_1 - \lambda_5)
(\lambda_5 - \lambda_4)(\lambda_4 - \lambda_1)},
\end{align*}
that is $c_1 = c_4$. By the same calculation using
 (\ref{ident2})--(\ref{ident8}), we obtain
\begin{align*}
 c_1 = c_4, \quad c_2 = - c_3, \quad c_5 = c_6, \quad c_1 = c_3, \quad 
c_2 = - c_4, \quad d_5 = d_6, \quad c_1 = - c_2, \\ 
d_2 = - \z d_6, \quad d_1 = - \z^3 c_5, \quad c_3 = - \z^4 d_4, \quad
c_4 = \z^4 d_3, \quad d_2 = d_3, \quad c_6 = \z^3 c_4.
\end{align*}
These equalities gives the ratios in the assertion.
\end{proof}
\begin{rem}
We have the following equalities;
\[
 J(ijklm) = - J(mlkji), \quad 
\Theta_{(a_1,a_2,a_3)}(\Omega(\eta))^5 = 
- \Theta_{(10 - a_1,10 - a_2,10 - a_3)}(\Omega(\eta))^5.
\]
\end{rem}
\vskip3mm
Let us denote $\x$ by $X$, and let $K_X$ be the canonical class of $X$.
\begin{cor}
We have an isomorphism of $\C$--algebras
\[
 \C[\Theta_{(a_1,a_2,a_3)}(\Omega(\eta))^5] \cong \oplus_{n = 0}^{\infty}
\mathrm{H}^0(X,\mathcal{O}_X(-nK_X)),
\]
where the left hand side is the $\C$--algebra of the functions on $\B$ 
generated by the twelve theta functions in Theorem \ref{maintheorem}. 
Especially the $\C$--vector space spanned by 
$\{ \Theta_{(a_1,a_2,a_3)}(\Omega(\eta))^5 \}$ coincides with 
$\mathrm{H}^0(X,\mathcal{O}_X(-K_X))$.
\end{cor}
\begin{proof}
The map $J$ is essentially anti-canonical map 
(see Section \ref{sectionconfig}). Hence it follows from Theorem
 \ref{maintheorem}.
\end{proof}
\begin{rem}
By the Riemann-Roch theorem, we obtain
\begin{align*}
\dim \mathrm{H}^0(X,\mathcal{O}_X(-n K_X)) &= 
\frac{1}{2}(-n K_X) \cdot (-n K_X -K_X) + 1 \\
&= \frac{5}{2} n(n+1) +1
\end{align*}
since $(-K_X) \cdot (-K_X) = 5$. So we have 
$\dim \mathrm{H}^0(X,\mathcal{O}_X(-K_X)) = 6$, and 
twelve $\Theta_{(a_1,a_2,a_3)}(\Omega(\eta))^5$ satisfy 6 linear
 equations. It is known that the image of $X$ in $\mathbb{P}^5$ 
by the anti-canonical map is determined by the system of quadratic 
equations (see \cite[Chapter 5]{F}).
\end{rem}
\vskip3mm
\noindent
{\bf ----- The graded ring of Automorphic forms----}
\par
Recall the automorphic factor $F_g(\eta)$ in Lemma \ref{lemcocycle}.
We consider the automorphic function $f(\eta)$ on $\B$ 
in the sense that we have
\begin{align} \label{modular}
f(g \eta) = F_g(\eta)^k f(\eta)\quad \text{for} \ g \in \Gamma(1-\z),
\end{align}
where $k$ is a non negative integer. Let us denotes the vector space of
holomorphic functions satisfying (\ref{modular}) by $A_k(F_g)$.
\vskip3mm
\begin{prop} \label{propmod}
Let $(a_1,a_2,a_3)$ and $(b_1,b_2,b_3)$ be the member of the system 
in (\ref{twelve1}), then it holds
\[
 \Theta_{(a_1,a_2,a_3)}(\Omega(\eta))^5 
\Theta_{(b_1,b_2,b_3)}(\Omega(\eta))^5 \in A_5(F_g).
\]
\end{prop}
\begin{proof}
By (\ref{newformula}) and Lemma \ref{lemmaphi'}, we have
\begin{equation} \label{eeee}
\begin{split}
&\Theta_{(a_1,a_2,a_3)}(\Omega(g \eta))^5 
\Theta_{(b_1,b_2,b_3)}(\Omega(g \eta))^5 \\
=&\kappa(g)^{10}\exp(2 \pi \sqrt{-1}\phi'_{(a_1,a_2,a_3)}(g))^{10}
F_g(\eta)^5 \Theta_{(a_1,a_2,a_3)}(\Omega(\eta))^5 
\Theta_{(b_1,b_2,b_3)}(\Omega(g \eta))^5
\end{split}
\end{equation}
for $g \in \Gamma(1-\z)$. We must show 
\begin{align} \label{kapp}
 [\kappa(g)\exp(2 \pi \sqrt{-1}\phi'_{(a_1,a_2,a_3)}(g))]^{10} = 1
\end{align}
for $g \in \Gamma(1-\z)$. Let $\eta_0$ be the point 
\[ (\text{the mirror of} \ h_{12}) \cap 
(\text{the mirror of} \ h_{34}) = [0:0:1] \in \B.
\] 
Then $\eta_0$ is fixed by $h_{12}$ and $h_{34}$. 
Moreover, $F_g(\eta) = \z$ for $h_{12}$ and $h_{34}$. So we have
\[
\Theta_{(1,1,1)}(\Omega(\eta_0))^{10} 
=[\kappa(g)\exp(2 \pi \sqrt{-1}\phi'_{(1,1,1)}(g))]^{10}
\Theta_{(1,1,1)}(\Omega(\eta_0))^{10} \quad (g = h_{12}, h_{34}) 
\]
by (\ref{eeee}). Since $\Theta_{(1,1,1)}(\Omega(\eta_0)) \ne 0$
(see Proposition \ref{evaluation}), we obtain (\ref{kapp}) 
for $h_{12}$ and $h_{34}$. By the same way, we see that (\ref{kapp})
 holds for any member $h_{ij}$ of the generator system of $\Gamma(1-\z)$. 
Hence it holds for all $g \in \Gamma(1-\z)$. 
\end{proof}
\begin{tm} \label{tmmodring}
(1). \ We have the isomorphism of the $\C$--algebras:
\begin{align*}
\oplus_{n = 0}^{\infty} A_{5n}(F_g) &= 
\C[\Theta_{(a_1,a_2,a_3)}(\Omega(\eta))^5 
\Theta_{(b_1,b_2,b_3)}(\Omega(\eta))^5] \\ &\cong 
\oplus_{n = 0}^{\infty} \mathrm{H}^0(X, \mathcal{O}_X(-2nK_X)). 
\end{align*}
(2). \ $A_n(F_g) = \{0\} \quad \text{for} \ n \in \mathbb{N},\  n 
\equiv 1,2,3,4 \ \text{mod} \ 5$. 
\end{tm}
\begin{proof}
By Proposition \ref{propmod}, a function $f \in A_5(F_g)$ defines the 
meromorphic function 
\[
 \frac{f(\eta)}{\Theta_{(1,1,1)}(\Omega(\eta))^{10}} 
\]
on $\B / \Gamma(1-\z)$. So, by Theorem \ref{maintheorem}, we have the
 isomorphism of $\C$--vector space:
\[
A_{5n}(F_g) \cong \mathrm{H}^0(X, \mathcal{O}_X(-2nK_X)) \quad \text{for} \ 
n \in \mathbb{N}.
\]
Hence we have the assertion (1). 
\par
Next let us recall that $X$ is the blow up of  $\mathbb{P}^2$ 
at 4 points. We denote this blow up by 
$\pi : X \rightarrow \mathbb{P}^2$. Then the Neron--Severi group 
$\mathrm{NS}(X)$ has the free generator $E_1, E_2, E_3, E_4$ and 
$\pi^*H$, where $\{ E_i \}$ are the exceptional curves with 
respect to $\pi$, and $H$ is a general line on $\mathbb{P}^2$. 
For $n \notin 5 \Z$, there is no 
divisor $D$ on $X$ such that $5 D = -2n K_X$ since 
$-K_X = 3 \pi^* H - E_1 - E_2 - E_3 - E_4$. This implies the assertion (2)
 since 
\[ 
A_n(F_g)^5 \subset A_{5n}(F_g) \cong 
\mathrm{H}^0(X, \mathcal{O}_X(-2nK_X)).
\]
\end{proof}
\vskip3mm
\noindent
{\bf ----- The Schwarz inverse for the Gauss HGDE 
$E_{2,1}(\frac{1}{5}, \frac{2}{5}, \frac{4}{5})$  ----}
\par
Let us consider the 1--dimensional disk
\[
 \mathbb{B}_1 = \{ \eta \in \B \ : \ \eta_1 = 0\},
\]
and the degenerate period map
\begin{align*}
 \Phi_{12} \ : L(12) \cong \pr \longrightarrow \mathbb{B}_1, \quad
t \mapsto [&0:\int_{\gamma_2} \omega:\int_{\gamma_3} \omega], \\
&\omega = x^{-\frac{4}{5}} (x-1)^{-\frac{2}{5}} (x-t)^{-\frac{2}{5}}dx,
\end{align*}
as in Section \ref{sectiondl} (the parameter $\lambda$ is specialized as
$(\lambda_1,\cdots,\lambda_5) = (0,0,1,t,\infty)$). Set
\[
 \Gamma(1-\z)_1 = 
\{ g \in \Gamma(1-\z) \ : \ g(\mathbb{B}_1) = \mathbb{B}_1\}.
\]
As we mentioned in Section \ref{sectiondl}, 
this is the triangle group $\Delta (5,5,5)$ up to the center.
Recall those are the Schwarz map and the monodromy group for Gauss 
hypergeometric differential equation 
$E_{2,1}(\frac{1}{5}, \frac{2}{5}, \frac{4}{5})$ (see
(\ref{Gauss(5,5,5)})). We have the explicit discription of the inverse:
\begin{tm} \label{Gaussinverse}
The map 
\[
\varTheta_{12} \ : \ \mathbb{B}_1 /  \Gamma(1-\z)_1 \longrightarrow
\pr, \quad \eta \mapsto 
[\Theta_{(1,1,1)}(\Omega(\eta))^5 : -\Theta_{(1,1,9)}(\Omega(\eta))^5]
\]
is an isomorphism, and this is the inverse map of the Schwarz map
\[
 \Phi_{12} \ : \pr \longrightarrow \mathbb{B}_1 / \Gamma(1-\z)_1, \quad
[1:t] \mapsto [0:\int_{\gamma_2} \omega:\int_{\gamma_3} \omega].
\]
\end{tm}
\begin{proof}
By Theorem \ref{maintheorem}, the restriction of meromorphic function
\[
 \frac{\Theta_{(1,1,9)}(\Omega(\eta))^5}{\Theta_{(1,1,1)}(\Omega(\eta))^5}
\]
on $L(12)$ is of order 1. In fact,
$L(12) \cap L(13) = L(12) \cap L(14) = L(12) \cap L(15) =  \phi$, so
the numerator vanishes at only $L(12) \cap L(35)$ with order 1,
the denominator vanishes at only $L(12) \cap L(45)$ with order 1, and
$L(12) \cap L(35) \ne L(12) \cap L(45)$ (see Section
 \ref{sectionconfig}). Hence the map $\varTheta_{12}$ is an isomorphism.
 Moreover, by Theorem \ref{maintheorem}, we have the equality
\[
 \frac{\Theta_{(1,1,9)}(\Omega(\eta))^5}{\Theta_{(1,1,1)}(\Omega(\eta))^5}
= - \frac{(\lambda_1 - \lambda_3)(\lambda_3 - \lambda_5)
(\lambda_5 - \lambda_2)(\lambda_2 - \lambda_4)(\lambda_4 - \lambda_1)}
{(\lambda_1 - \lambda_3)(\lambda_3 - \lambda_2)
(\lambda_2 - \lambda_4)(\lambda_4 - \lambda_5)(\lambda_5 - \lambda_1)},
\]
on $\B / \Gamma(1-\z) \cong \x$, and this induces the equality
\[
 \frac{\Theta_{(1,1,9)}(\Omega(\eta))^5}{\Theta_{(1,1,1)}(\Omega(\eta))^5}
= - \frac{(\lambda_3 - \lambda_5)(\lambda_4 - \lambda_1)}
{(\lambda_3 - \lambda_1)(\lambda_4 - \lambda_5)}
\]
on $L(12)$. Putting $(\lambda_1,\lambda_3,\lambda_4,\lambda_5) = 
(0,1,t,\infty)$, we obtain
\[
 \frac{(\lambda_3 - \lambda_5)(\lambda_4 - \lambda_1)}
{(\lambda_3 - \lambda_1)(\lambda_4 - \lambda_5)} = t.
\]
\end{proof}
\vskip3mm
Let us consider a holomorphic function $f$ on $\mathbb{B}_1$ satisfying
the condition:
\[
 f(g \eta) = F_g(\eta)^k f(\eta)\quad \text{for} \ g \in \Gamma(1-\z)_1,
\]
and we denote the $\C$--vector space of such functions by $M_k(F_g)$.
\vskip3mm
\begin{cor}
(1). \ We have an isomorphism of $\C$--algebras:
\begin{align*}
&\oplus_{n = 0}^{\infty} M_{5n}(F_g) \\ 
= &\C[\Theta_{(1,1,1)}(\Omega(\eta))^{10},\ 
\Theta_{(1,1,1)}(\Omega(\eta))^5 \Theta_{(1,1,9)}(\Omega(\eta))^5, \
\Theta_{(1,1,9)}(\Omega(\eta))^{10}] \\ 
\cong &\C[x_0^2, \ x_0x_1, \ x_1^2] \\
\cong &\oplus_{n = 0}^{\infty} 
\mathrm{H}^0(\pr, \mathcal{O}_{\pr}(-nK_{\pr})), 
\end{align*}
where $[x_0:x_1]$ is homogeneous coordinates of $\pr$. \\
(2). \ $M_n(F_g) = \{0\} \quad \text{for} \ n \in \mathbb{N},\  n 
\equiv 1,2,3,4 \ \text{mod} \ 5$. 
\end{cor}
\begin{proof}
The assertion (1) is a direct consequence of Corollary \ref{tmmodring}
 and Theorem \ref{Gaussinverse}. The assertion (2) is obtained by the
 same argument as the proof of Theorem \ref{tmmodring}.
\end{proof}
\vskip3mm
\noindent
{\bf Acknowledgments.} 
I express my sincere thanks to Professor Hironori Shiga for
advices during the prepartion of this paper.

\end{document}